\documentclass[thmsa,final,french,11pt]{article}

\pdfoutput=1
\usepackage{ifpdf}
\expandafter\ifx\csname ifpdf\endcsname\relax
 \csname newif\expandafter\endcsname\csname ifpdf\endcsname
 \else

\usepackage{eurosym}
\usepackage{amsfonts}
\usepackage{amsmath}
\usepackage{amssymb}
\usepackage{mathrsfs} 
\usepackage{enumerate}
\usepackage{graphicx}

\usepackage{ntheorem}

\newtheorem{Th}{Theorem}

\newtheorem{Prop}[Th]{Proposition}
\newtheorem{Le}[Th]{Lemma}
\newtheorem{Cor}[Th]{Corollary}

\newcommand{\WW}{\mathbf{W}}

\newcommand{\pth}[1]{\ensuremath{\left(#1\right)}}
\newcommand{\floor}[1]{\ensuremath{\left\lfloor#1\right\rfloor}}

\newcommand{\diff}{{\stackrel{\neq}{\ ..}}}
\newcommand{\size}[1]{\ensuremath{\card\pth{#1}}}

\newcommand{\length}[1]{\ensuremath{\ell\pth{#1}}}

\newcommand{\pbrcx}[1]{\ensuremath{\left[ {#1} \right]}}

\DeclareMathOperator\card{Card}

\DeclareMathOperator{\Probability}{\mathbb{P}}

\DeclareMathOperator{\Areaa}{\mathcal A}
\DeclareMathOperator{\Expected}{\mathbb{E}}

\newcommand{\Ex}[1]{\Expected\pbrcx{#1}}
\newcommand{\VVV}[1]{\mathbb{V}\left[#1 \right]}
\newcommand{\Prob}[1]{\Probability\pbrcx{#1}}

\newcommand{\Area}[1]{\Areaa\pth{#1}}
\newcommand{\eangle}[1]{|\widehat{#1}|}
\newcommand{\expo}[1]{e^{#1}}
\newcommand{\R}{{\mathbb{R}}}
\newcommand{\s}{{\mathbb{S}}}

\newcommand{\N}{{\mathbb{N}}}
\newcommand{\Z}{{\mathbb{Z}}}
\newcommand{\setpath}{{\mathcal{P}_{s,t}(X)}}
\newcommand{\indicator}[1]{\mathbf{1}_{[#1]}}

\newcommand{\midd}{\;\rule[-4mm]{0.5mm}{10mm}\;}

\newcommand{\DT}{\operatorname{Del}}

\newcommand{\grid}{{\mathbf G}}
\newcommand{\animal}{{\mathbf A}}

\newcommand{\betamatrix}{\mbox{\footnotesize$\begin{pmatrix}
1 & 1 & 1\\ \cos\beta_1 & \cos\beta_2 & \cos\beta_3\\ \sin\beta_1 & \sin\beta_2 & \sin\beta_3 
\end{pmatrix}$}}
\newcommand{\smallasin}{\mbox{\tiny{$\arcsin\,$}}}

\def\WW_#1{\boldsymbol{W}\!_#1}
\numberwithin{equation}{section}

\newenvironment{prooft}{\noindent {\bf Proof}} {\hfill $\square$  \noindent}

\usepackage{xcolor}

\usepackage[french]{varioref}
\begin{document}

\title{Stretch Factor of Long Paths in a planar Poisson-Delaunay Triangulation\thanks{Partially supported by ANR blanc PRESAGE
    (ANR-11-BS02-003).}}
\author{Nicolas CHENAVIER\footnote{Universit\'e Littoral C\^ote d'Opale, EA 2797, LMPA, 50 rue Ferdinand Buisson, F-62228 Calais, France.}, Olivier DEVILLERS\footnote{Inria, Centre de Recherche Nancy-Grand Est, France.}\footnote{CNRS, Loria, France.}\footnote{Université de Lorraine, France.}}
\maketitle

\begin{abstract}
Let $X:=X_n\cup\{(0,0),(1,0)\}$, where $X_n$ is a planar
  Poisson point process of intensity $n$.
We provide a first non-trivial lower bound for
the distance between the expected length of the shortest path between
$(0,0)$ and $(1,0)$ in  the Delaunay triangulation associated with $X$
when the  intensity of $X_n$ goes to infinity. 
Experimental values indicate that the correct value is about 1.04.
 We also prove that the expected number of Delaunay edges crossed
  by the line segment $[(0,0),(1,0)]$ is equivalent  to
  $2.16\sqrt{n}$ and that the expected length of a particular path converges
  to 1.18 giving an upper bound on the stretch factor.
\end{abstract}

%\strut


%\strut


\section{Introduction}
Let $\chi$ be a locally finite subset in $\R^d$, endowed with its Euclidean norm $\|\cdot \|$, such that each subset
of size $n<d+1$ are affinely independent and no $d+2$ points lie on a
sphere. If $d+1$ points $x_1,\ldots, x_{d+1}$ of $\chi$ lie on a ball
that contains no point of $\chi$ in its interior, we define an edge
between $x_i$ and $x_j$ for each $1\leq i\neq j\leq d+1$. The set of
these edges is denoted by $\DT(\chi)$ and the graph $(\chi,\DT(\chi))$
is the so-called Delaunay graph associated with $\chi$~\cite[p.~478]{schneider08}. 
Delaunay triangulation is a very popular structure in
  computational geometry~\cite{aurenhammer2013voronoi}
and is extensively used in many areas such as surface reconstruction~\cite{cazals2006delaunay}
or mesh generation~\cite{cheng2012delaunay}. 

In this paper, we investigate {several}  paths in $\DT(\chi)$. By a path
$P=P(s,t)$ between two {points} $s,t\in \chi$, we mean a sequence of
segments $[Z_0, Z_1], [Z_1,Z_2], \ldots, [Z_{k-1},Z_k]$, {such that} $Z_0=s$, $Z_k=t$. {In particular, we say that $P$ is a path in $\DT(\chi)$ if it is a path such that each segment $[Z_i, Z_{i+1}]$ is an edge in $\DT(\chi)$}. {The investigation of paths} is related to walking
strategies which are commonly used to perform point location in planar 
triangulations~\cite{devillers02} or routing in geometric
networks~\cite{bose04}. One of the classical works concerns the so-called straight walk which deals with the set of triangles cut by the line segment $[s,t]$. In this context,  Devroye, Lemaire and
 Moreau~\cite{devroye04} consider $n$ points evenly distributed in a
 convex domain and prove that the expected number of Delaunay edges 
crossed by a line segment of length $L$  is $O\pth{L\sqrt{n}}$. This
result is improved by Bose and Devroye~\cite{bose06} who show that the
complexity equals $\Theta(\sqrt{n})$, 
even if
the line segment depends on the
distribution. 

\begin{figure}[t]
     \begin{center}
         \includegraphics[page=1]{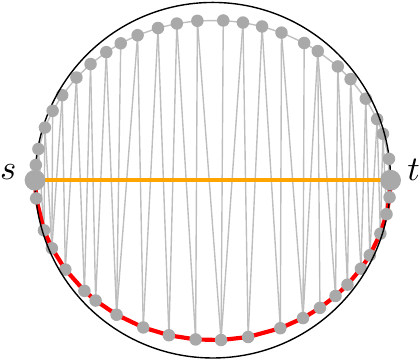}
     \end{center}
     \caption{
        The shortest path in the Delaunay triangulation (red) can be
        about $\frac{\pi}{2}$ larger than the line segment (yellow).
        \label{fig:pi_over_two}
     }
 \end{figure}

Another classical problem in walking strategies is the investigation
of the stretch factor associated with two nodes $s,t\in \chi$ in $\DT(\chi)$. This quantity is defined as
$\frac{l(SP_{\chi})}{\|s-t\|}$, where $l(SP_{\chi})$ is the length of
the shortest path between $s$ and $t$. Many
upper bounds were established for the stretch factor in the context of
finite sets $\chi$, e.g.~\cite{dfs-dgaag-90, kg-dtcac-89}.
The best upper bound established until now
for deterministic
finite sets $\chi$ is due to  Xia~\cite{x-sfdtl-13} who
proves that the stretch factor is lower than 1.998.   
For the lower bound, Xia and Zhang~\cite{xz-ttbsf-11} find a
configuration of points $\chi$ such that the stretch factor
is greater than 
1.5932 (e.g. Figure~\ref{fig:pi_over_two} provides a
configuration where the stretch factor is close to
$\pi/2\simeq1.5708$).  

In this paper, we focus on a probabilistic version  
of the problem
 by taking a slight modification of the underlying point process.
More precisely, we consider a homogeneous Poisson point process $X_n$ of
 intensity $n$ in all the plane $\R^2$.
For such an infinite point process, studying the maximum of the
stretch over any points $s,t\in X_n$ has no real sense. Indeed, for
any bounded set $B\subset\R^2$, there exists a configuration of points
in $X_n\cap B$ close to the one depicted in
Figure~\ref{fig:pi_over_two} with non-zero probability. In particular,
such a configuration occurs almost surely somewhere in the plane.  

\paragraph{Contributions}
The main focus of our paper is
to provide bounds for the expectation of the stretch factor
between two fixed points $s,t\in\R^2$ in the Delaunay triangulation
$\DT(X)$, where $X:=X_n\cup\{s,t\}$.  The difficulty to obtain a lower bound for $\Ex{\length{SP_{X}}}$ comes
from the fact that we absolutely do not know where the shortest path $SP_{X}$
is. We take the challenge by establishing our first main theorem.

\begin{Th} 
\label{th:lower}
Let $X:=X_{n}\cup\{s,t\}$, where $X_n$ is a Poisson point process of intensity $n$ and $s,t\in\R^2$. Then 
\[ \Prob{\length{SP_{X}} \le (1+2.47\cdot 10^{-11}) \|s-t\|  }  =  O\pth{n^{-\tfrac{1}{2}}}.\]
\end{Th}
As a consequence, we easily {deduce} the following result:
\begin{Cor}
Let $X:=X_{n}\cup\{s,t\}$, where $X_n$ is a Poisson point process of intensity $n$ and $s,t\in\R^2$. Then 
\[\liminf_{n\rightarrow\infty}\Ex{\frac{\length{SP_{X}}}{\|s-t\|}}\geq  1+2.47\cdot 10^{-11}.\]
\label{cor:lower}
\end{Cor}
We think that our results provide the first non-trivial lower bound
(i.e. greater than 1) for the stretch factor when the intensity of the
underlying Poisson point process goes to infinity. However, our lower
bound is far from optimal since simulations suggest that
$\lim_{n\rightarrow\infty}\Ex{\frac{\length{SP_{X}}}{\|s-t\|}}\simeq
1.04$. We notice that our result is closely related to a theorem
recently proved by Hirsch, Neuhaüser, and Schmidt~\cite[Theorem~26]{HNS}.
 Indeed, they show
that $\inf_{n\geq 1}\Ex{\frac{l(SP_{X})}{\|s-t\|}}>1$. However, 
%they do not  
their technique cannot 
provide  explicit lower bound for the stretch factor. In the
following proposition, we also give an upper bound.
\begin{Th}\label{prop:upper-bound}
Let $X:=X_{n}\cup\{s,t\}$, where $X_n$ is a Poisson point process of intensity $n$ and $s,t\in\R^2$. Then 
\[
\limsup_{n\rightarrow\infty}\Ex{\frac{\length{SP_{X}}}{\|s - t\|}}
\leq   \frac{35}{3\pi^2} \simeq 1.182.
\]
\end{Th}
The upper bound we considered above is established by bounding the
length of a particular path in the Delaunay triangulation. In
particular, our theorem improves a result due to Baccelli et
al.~\cite{BTZ}. Indeed, by considering another particular path, they
prove that the {expectation of the } stretch factor is lower than 
$\frac{4}{\pi}\simeq 1.27$.

\paragraph{Outline}
In Section~\ref{s:preliminaries}, we begin with some preliminaries by introducing notation and tools of stochastic and integral geometry. As a first result, we use these tools in the next section to obtain a  tight evaluation of the size of the straight walk. In particular, this
evaluation gives an explicit value for the constant appearing in the work by
Bose and Devroye~\cite{bose06}
 which was hidden in the asymptotic complexities. In
 Section~\ref{s:upperBound}, we provide estimates for the length of a
 particular path. These estimates directly imply Theorem~\ref{prop:upper-bound}. Section~\ref{s:lowerBound} constitutes the
 main part of our paper and deals with the lower bound for the
 shortest path. Our main idea is to discretize the plane into pixels
 and to consider the so-called lattice animals. We derive Theorem~\ref{th:lower} by investigating the size of these lattice animals and by adapting tools of percolation theory. In Section~\ref{s:simulation}, we give experimental values {for} various quantities. These simulations confirm the results for the size and the
length of the straight walk and  suggest that the correct value
of the expected stretch factor is 1.04.  In Appendix, we give auxiliary results which are used throughout the paper.

\section{Preliminaries}
\label{s:preliminaries}

\paragraph*{Notation}

Let $s$ and $t$ be two fixed points in $\R^2$.
Let $X_n$ be a homogenous Poisson point process of intensity $n$ in $\R^2$
and $X:=X_n\cup\{s,t\}$. We denote by $\DT(X)$ the Delaunay triangulation associated with $X$.  We give below several notation which will be used throughout the
 paper. 
\begin{itemize}
\item For any point $p\in\R^2$, we write $p=(x_p,y_p)$.
\item For any segment $e\subset\R^2$, we denote by $h(e)$ and
  $\eangle{e}$ the length of the horizontal projection of $e$ and the
  absolute value of the angle with the $x$-axis respectively.  
\item  For any $k$-tuple of points $p_1, \ldots, p_k\in\R^2$, we write 
 $p_{1:k}:=(p_1,\ldots, p_k)$. When $p_1, \ldots, p_k\in\R^2$ are
 pairwise distinct, we write the $k$-tuple of points as $p_{1\diff k}$. 
Such a notation will be used in the summation index. Moreover,
 with a slight abuse of notation, we write $\{p_{1:k}\}:=\{p_1,\ldots, p_k\}$. 
 \item For each 3-tuple of points $p_{1\diff 3}\in (\R^2)^3$ which do not belong to the same line, we denote by $\Delta(p_{1: 3})$
   and $B(p_{1: 3})$ the triangle spanned by $p_{1: 3}$ and
   the (open) circumdisk associated with $p_{1: 3}$ respectively. 
  We also denote by $R(p_{1: 3})$ the radius of $B(p_{1: 3})$.
 \item For each $z\in\R^2$ and $r\geq 0$, let $B(z,r)$ be the disk centered at $z$ with radius $r$. 
 \item Let $\s\subset \R^2$ be the unit circle endowed with its 
  uniform distribution  $\sigma$ such that $\sigma(\s)=2\pi$.
  \item For any Borel subset $B\subset\R^2$, let $\Area{B}$ be the area of $B$.
 In particular, for each $u_{1:3}\in \s^3$, we have
\[
\Area{\Delta(u_{1:3})} = \tfrac{1}{2} \left| \det\betamatrix \right| ,
\]
where $\beta_i\in[0,2\pi)$ is the angle between $u_i$ and $(1,0)$, with $1\leq i\leq 3$.
\end{itemize}

\paragraph*{Paths}

\begin{figure}[t]
     \begin{center}
         \includegraphics[page=2,width=\textwidth]{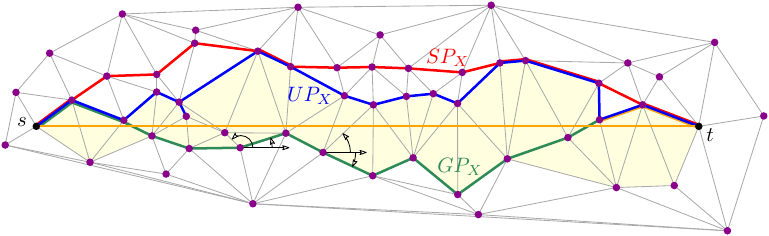}
     \end{center}
     \caption{
        The paths $SP_X$, $UP_X$, and $GP_X$  
         (the path $GP_X$  is defined in Section~\ref{s:simulation}).
        \label{fig:paths}
     }
 \end{figure}
 Given $s,t\in\R^2$, we denote by  $\setpath$ the family of paths in $\DT(X)$ starting from $s$ and going to $t$, with $X:=X_n\cup\{s,t\}$. For each path $P\in\setpath$, we denote by
$\length{P}$ the Euclidean length of $P$
and $\size{P}$ its number of edges. We introduce {two} types of paths in $\setpath$, namely 
$SP_X$ and $UP_X$ {as follows:} 

 \label{def:SP}
 \textit{Shortest Path $SP_X$:} this path is the one which minimizes the length between $s$ and $t$ 
 in the Delaunay triangulation $\DT(X)$.

\textit{Upper Path $UP_X$:} this path is defined as the set of all edges 
in $\R\times\R_+$ which belong to Delaunay triangles that intersects $[s,t]$.
Some of these edges may be traversed in both ways (e.g. this is the case for the edges
incident to the fifth vertex of the blue path in
Figure~\ref{fig:paths}). 

These {two} paths are depicted in 
Figure~\ref{fig:paths} and described below.
{A third path, referred to as the greedy path $GP_X$, is also depicted in  Figure~\ref{fig:paths} but introduced in Section \ref{s:simulation}. } 
%this is the path defined  as the piece of the boundary of
%  $\bigcup_{i=0}^k T_i$ above the line $[s,t]$,
%where  $(T_i)_{0\leq i\leq k}$  is the family
%of triangles in $\DT(X)$  whose interior  intersects $[s,t]$ and ordered from the
%left to the right, with $k\geq 1$. 

%%%%%%%%%% MOVED TO EXPERIMENTS
% \textit{Greedy constructed Path $GP_X$:} we define such a path, starting
% in $s$, by induction. 
% Let $(T_i)_{0\leq i\leq k}$  be the family
% of triangles in $\DT(X)$  whose interior  intersects $[s,t]$ and ordered from the
% left to the right.
% Let $w$ be the last vertex inserted in the path and let $i$ be the greatest value such that  $w\in T_i$. 
% Then the edge of $T_i$ incident to $w$ which minimizes the angle with respect to 
% the $x$-axis is added to the path.

\paragraph*{Slivnyak-Mecke and Blaschke-Petkantschin Formulas }

Throughout the paper, we will extensively use two classical formulas
of stochastic and integral geometry. As a warm-up, we first provide a new proof of a well-known result (usually established with the Euler relation) to introduce the computation technique with an easy pedagogical case.

\begin{Prop}\label{prop:origin}
Let $X_n$ be a Poisson point process of intensity $n$ and let $N_0$ be the number of triangles in $\DT(X_n)$  whose circumdisk contains the origin. Then $\Ex{N_0}=4$. 
\end{Prop}
\begin{prooft}
First we notice that
\begin{eqnarray*}
\Ex{  N_0}  &=& 
{\tfrac{1}{3!}}
\Ex{\sum_{{p_{1\diff 3}\in X_n^3}}
       \indicator{\Delta(p_{1: 3})\in \DT(X_n)}
       \indicator{O\in B(p_{1: 3})} }. 
\end{eqnarray*}
 A first tool allows us to re-write the right-hand side as an
 integral. More precisely, according to the  Slivnyak-Mecke formula 
(e.g.~\cite[Theorem~3.3.5]{schneider08}), we have 
\begin{eqnarray*}
\Ex{N_0}  &=& 
\tfrac{n^3}{3!}\int_{(\R^2)^3} \Prob{\Delta(p_{1:3})\in \DT(X_n\cup\{p_{1:3}\})}
        \indicator{O\in B(p_{1:3})}      dp_{1:3}
\\ &=&
\tfrac{n^3}{6}  \int_{(\R^2)^3} \expo{-n\Area{B(p_{1:3})}}   \indicator{O\in B(p_{1:3})}  dp_{1:3}
\end{eqnarray*}
since $X_n$ is a Poisson point process. A second formula transforms
the integral  over $(\R^2)^3$ as an integral over $\R^2\times \R_+\times
\s^3$ by associating with each $p_{1:3}\in (\R^2)^3$ the circumcenter, the
circumradius and the angles of $p_i$, $1\leq i\leq 3$,
respectively. More precisely, from the Blaschke-Petkantschin formula
(e.g.~\cite[Theorem~7.3.1 ]{schneider08}), we have
\begin{eqnarray*}
\Ex{ N_0}  &=& 
\tfrac{n^3}{6} \int_{\R_+} \int_{\R^2} \int_{\s^3 }
        \expo{-n\Area{B(z,r)}}
       \indicator{O\in B(z,r)} 
\,\cdot\,
       r^3
          2\Area{\Delta(u_{1:3})}
        \sigma(du_{1:3}) dz dr. 
\end{eqnarray*}
Integrating over $z$ by noting that $O\in B(z,r)\Longleftrightarrow
z\in B(O,r)$, we get
\begin{equation*}
\Ex{ N_0} = \tfrac{\pi n^3}{{6}}\int_{\R_+} {\expo{-n\pi r^2}}  r^5dr \times \int_{\s^3} {2}\Area{\Delta(u_{1:3})}
        \sigma(du_{1:3})
\end{equation*} 
since $\int_{B(O,r)}dz=\pi r^2$ for each $r\geq 0$. It follows that

\begin{align*}
\Ex{ N_0} 
& = \frac{\pi n^3}{6}\int_{\R_+}{\expo{-n\pi r^2}}  r^5dr \times \int_{[0,2\pi)^3} \left| \det\betamatrix   \right| d\beta_{1:3} \\
& = \frac{\pi n^3}{6}\times 
\frac{1}{\pi^3n^3}\times 
24\pi^2= 4.
\end{align*}
\end{prooft}

A table of integrals is provided in Appendix~\ref{ap:integrals}. As a second easy warm up, we deal with the sum of the lengths of the
Delaunay edges with respect to a typical vertex. In the following proposition, we only give an upper-bound for such a length since it is enough for us in the sequel.

\begin{Prop}\label{prop:origin-length}
Let $X_n$ be a Poisson point process of intensity $n$ and let $L_0$ be
the sum of the lengths of the edges with vertex $O$ in 
$\DT(X_n\cup\{O\})$, where $O\in\R^2$ denotes the origin. Then
$\Ex{L_0} =  c\cdot n^{-\tfrac{1}{2}}$ for some constant $c$. 
\end{Prop}
In the above proposition, we do not make explicit the constant $c$. However, experimental values suggest that $c\simeq 6.8$. 

\begin{prooft}
The sum of the length of the edges $L_0$ can be expressed as
\[
\Ex{  L_0}  =
\tfrac{1}{2}\tfrac{1}{2}
\Ex{\sum_{{p_{1\diff 2}\in X_n^2}}
       \indicator{\Delta(p_1p_2O)\in \DT(X_n\cup\{O\})}
        \pth{ \|p_1\|+\|p_2\|}}. 
\]
The first factor $\tfrac{1}{2}$ comes from the fact that each triangle  
is  counted twice (once clockwise and once counterclockwise). The second factor $\tfrac{1}{2}$ results from the fact that each edge is also obtained twice (once for each incident triangle). It follows that 
\begin{equation}
\label{eq:defL0}
\Ex{L_0}  = 
\tfrac{n^2}{4}\int_{(\R^2)^2} \Prob{\scriptstyle\Delta(p_{1}p_2O)\in \DT(X_n\cup\{p_{1},p_2,O\})}
            \pth{ \|p_1\|+\|p_2\|}   dp_1dp_2. 
\end{equation}
To make explicit the right-hand side, we give an analogous 
version of the Blaschke- Petkanschin type change of variables in which
one of the vertices is held fixed. We proceed in the same spirit 
as in Schneider and Weil~\cite[proof of Theorem~7.3.2]{schneider08}. 
More precisely, let $\phi$ be the function\\
\begin{minipage}{0.5\textwidth}
\[\begin{split} \phi: \;  &   \R_+\times [0,2\pi)^3\longrightarrow \R^2\times \R^2\\
& (r,\alpha, \beta_1,\beta_2) \longmapsto (p_1,p_2),
\end{split}
\] where, for each $i=1,2$, we let
\[p_i=r(\cos\alpha + \cos\beta_i, \sin\alpha+\sin\beta_i).\]
\end{minipage}\hfill\begin{minipage}{0.4\textwidth}
 \includegraphics[page=3,width=\textwidth]{Figures}
\end{minipage}\medskip\\
Provided that $O$, $p_1$ and $p_2$ do not belong to the same line, we 
notice that $r$ and $z:=r(\cos\alpha, \sin\alpha)$ are the
circumradius and the circumcenter of the triangle $\Delta(O,p_1,p_2)$
respectively. The terms $\beta_1$ and $\beta_2$ represent the angles
of the vectors 
$\overrightarrow{zp_1}$ and $\overrightarrow{zp_2}$
 with respect to the $x$-axis. Up to a null set, the function $\phi$ defines a $C^1$ diffeomorphism with Jacobian $J_\phi(r,\alpha,\beta_1,\beta_2):=r^3\cdot D(\alpha,\beta_1,\beta_2)$, where
\[D(\alpha,\beta_1,\beta_2):= \begin{vmatrix}
\cos\alpha+\cos\beta_1 & -\sin\alpha & -\sin\beta_1 & 0\\
\sin\alpha+\sin\beta_1 & \cos\alpha & \cos\beta_1 & 0\\
\cos\alpha + \cos\beta_2 & -\sin\alpha & 0 & -\sin\beta_2\\
\sin\alpha+\sin\beta_2 & \cos\alpha & 0 & \cos\beta_2
\end{vmatrix}  \]
Since 
%$\|p_i\|=\od{r}\sqrt{2}\sqrt{1+\cos(\alpha-\beta_i)}$ 
$\|p_i\|=2r\left|\cos\frac{\beta_i-\alpha}{2}\right|$
with $(p_1,p_2)=\phi(r,\alpha,\beta_1,\beta_2)$, it follows from  \eqref{eq:defL0} that
\begin{multline*}
\Ex{L_0} = \tfrac{n^2}{4}\int_{\R_+}\int_{[0,2\pi)^3}e^{-n\pi r^2}
2r \pth{ \left|\cos\frac{\beta_1-\alpha}{2}\right| + \left|   \cos\frac{\beta_2-\alpha}{2}\right| }\\
\times r^3\cdot \left|D(\alpha,\beta_1,\beta_2)   \right|d\alpha d\beta_1d\beta_2dr\\
=\tfrac{n^2}{4}\cdot \frac{3}{8\pi^2n^2\sqrt{n}}\;\;\;
            2\!\!\int_{[0,2\pi)^3} \!\!\!\mbox{\footnotesize$ \pth{ \left|\cos\!\frac{\beta_1-\alpha}{2}\right| + \left|   \cos\!\frac{\beta_2-\alpha}{2}\right| }
            \left|D(\alpha,\beta_1,\beta_2)   \right|d\alpha d\beta_1d\beta_2dr$}
\end{multline*}
It follows that $\Ex{L_0} = c\cdot n^{-\tfrac{1}{2}}$ since the integral over $[0,2\pi)^3$ is finite. 
\end{prooft}

\section{Size of Straight Walk\label{s:straight}}
\label{s:straightwalk}
In this section and in the next one, we restrict our interest to
triangles which are cut by the line segment $[s,t]$. We first
investigate the mean number of edges for the straight walk.

\begin{Prop}\label{lem:straight}
Let $X:=X_{n}\cup\{s,t\}$, where $X_n$ is a Poisson point process of intensity $n$ and $s,t\in\R^2$.  Let
$N_{s,t}$ be the number of edges of triangles in $\DT(X)$ intersecting $[s,t]$. Then
\[
\Ex{\frac{N_{s,t} }{\|s - t\|}} 
=\frac{64}{3\pi^2} \sqrt{n} \;+O(1) \simeq 2.1615 \sqrt{n} .
\]
\end{Prop}
Even if  the previous proposition is not closely related to the
shortest path, we think that our result could help us for a better
understanding of $SP_X$. In particular, the tools appearing in the
proof of Proposition~\ref{lem:straight} will be useful in the
following sections.

\begin{prooft}
Without loss of generality, we assume that $s=(0,0)$ and
$t=(1,0)$. First we notice that almost surely
\begin{equation}
\label{e:Nst}
N_{s,t} = \tfrac{1}{2} \sum_{p_{1\diff 3}\in X_n^3} \indicator{\Delta(p_{1:3})\in \DT(X)}
       \indicator{p_{1:3}\in E^+} + \tfrac{1}{2}
        \sum_{p_{1\diff 3}\in X_n^3} \indicator{\Delta(p_{1:3})\in \DT(X)}
       \indicator{p_{1:3}\in E^-} + 1,\end{equation}
       where 
  \begin{align*}
&E^+:=\{p_{1:3}\in (\R^2)^3: \Delta(p_{1:3})\cap [s,t]\neq \varnothing, B(p_{1:3})\cap\{s,t\}=\varnothing, y_{p_1}, y_{p_2}> 0, y_{p_3}<0\},
\\ &E^-:=\{p_{1:3}\in (\R^2)^3: \Delta(p_{1:3})\cap [s,t]\neq \varnothing, B(p_{1:3})\cap\{s,t\}=\varnothing, y_{p_1}, y_{p_2}< 0, y_{p_3}>0\}. 
\end{align*}
In the same spirit as above, we also define a third set which will be useful in the sequel as 
\[E'^+:=\{p_{1:3}\in (\R^2)^3: \Delta(p_{1:3})\cap \mbox{ line}(s,t)\neq \varnothing, y_{p_1}, y_{p_2}> 0, y_{p_3}<0\}.\]
Notice that the factor $\tfrac{1}{2}$ appearing in \eqref{e:Nst} comes from the fact that each triangle of $E$ is obtained
 twice, with counterclockwise and clockwise orientations. 
Indeed, the triangles crossed by $[s,t]$ have vertices in $X_n$ except the first and the last
ones which are incident to $s$ and $t$ respectively. Besides, the number of edges crossed
by $[s,t]$ is one less than the number of triangles {intersecting $[s,t]$}, which explains the term +1 in Equation~(\ref{e:Nst}).

Provided that $p_{1:3}\in E^+$, we notice that $\Delta(p_{1:3})$ is a triangle in $\DT(X)$ if and only if it is a triangle in $\DT(X_n)$. This implies that
\[ \Ex{\sum_{p_{1\diff 3}\in X_n^3} \indicator{\Delta(p_{1:3})\in \DT(X)}  \indicator{p_{1:3}\in E^+}}
= \Ex{\sum_{p_{1\diff 3}\in X_n^3} \indicator{\Delta(p_{1:3})\in    \DT(X_n)}  \indicator{p_{1:3}\in E^+}}. \]
Besides, thanks to the Slivnyak-Mecke formula, we have
\begin{equation*}
\Ex{\sum_{p_{1\diff 3}\in X_n^3} \indicator{\Delta(p_{1:3})\in \DT(X_n)}
       \indicator{p_{1:3}\in E^+}} = n^3 \int_{(\R^2)^3} 
%    \Prob{\sharp X \cap B(p_{1:3})=0}\indicator{p_{1:3}\in E^+}
    \Prob{{ X \cap B(p_{1:3})=\varnothing}}\indicator{p_{1:3}\in E^+}
      dp_{1:3}.
\end{equation*}
When $p_{1:3}\in E^+$, the circumdisk $B(p_{1:3})$ intersects $[s,t]$
and does not contain $s$ {or} $t$. In particular, we have
$z(p_{1:3})\in [0,1]\times[-R(p_{1:3}),R(p_{1:3})]$. 
It  follows from the Blaschke-Petkantschin formula that 
\begin{multline}
\label{e:majE}
 \Ex{\sum_{p_{1\diff 3}\in X_n^3} \indicator{\Delta(p_{1:3})\in \DT(X_n)}
       \indicator{p_{1:3}\in E^+}}  \\
=  n^3 \int_{0}^\infty \int_0^1\int_{-r}^r \int_{\s^3 } \expo{-n\pi r^2}
       \indicator{z+ru_{1:3}\in E^+}
       \,\cdot\,
         r^3  2\Area{\Delta(u_{1:3})}
        \sigma(du_{1:3}) dy_zdx_z dr. 
\end{multline}
\begin{figure}[t]
     \begin{center}
         \includegraphics[page=4,width=0.8\textwidth]{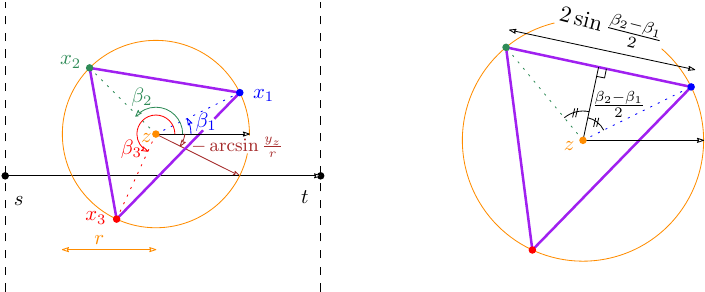}
     \end{center}
     \caption{
       Notation for the proofs of Propositions~\ref{lem:straight} and~\ref{lem:upper-bound}.
        \label{fig:notation}
     }
 \end{figure}
Let $r>0$ be fixed. First, we make explicit the integral over $y_z$
  and $u_{1:3}$ in the above equation when we replace the set $E^+$ by $E'^+$. By taking the
change of variables $h=\frac{y_z}{r}$ (see Figure~\ref{fig:notation}),
we have  
\begin{multline*}
\begin{split}
& \int_{-r}^r    \int_{\s^3} 
      2\Area{\Delta(u_{1:3})}\indicator{z+ru_{1:3}\in E'^+}\sigma(du_{1:3})dy_z   
\\ 
       & =  r \!\!
 \int_{-1}^1 
      \int_{\pi+\smallasin h}^{2\pi-\smallasin h}\!\!
     \int_{-\smallasin h}^{\pi+\smallasin h}\!\!
      \int_{-\smallasin h}^{\pi+\smallasin h}
\!\!   \left | \det \betamatrix \! \right|
      d\beta_1  d\beta_2  d\beta_3  dh
%\\   &
    =  \tfrac{512}{9}r .\hspace*{-8mm}\raisebox{2mm}{}
       \end{split}
\end{multline*}
%This implies that
We can evaluate the right-hand side of Equation~(\ref{e:majE}) {by substituting 
  $E^+$ by $E'^+$} as follows: 
\begin{multline}
\label{e:majE2}
 n^3 \int_{0}^\infty \int_0^1\int_{-r}^r \int_{\s^3 } \expo{-n\pi r^2}
       \indicator{z+ru_{1:3}\in E'^+}
       \,\cdot\,
         r^3  2\Area{\Delta(u_{1:3})}
        \sigma(du_{1:3}) dy_zdx_z dr\\
         =   \frac{512n^3}{9} \int_0^\infty    \int_0^1      \expo{-n\pi r^2}         r^4 dx_z dr =
\frac{512 n^3}{9} \cdot   \frac{3}{8\pi^2 n^2\sqrt{n}}= 
\frac{64}{3\pi^2} \sqrt{n}. 
\end{multline}

\noindent
\begin{minipage}{0.7\textwidth} 
To prove that the right-hand side of \eqref{e:majE} has the same order
when we replace $E^+$ by $E'^+$, we define:
\[
{ B^{\cup}(s,t,r):= \pth{B(s,r)\cap\{z; x_z\ge 0\} } \cup \pth{B(t,r)\cap\{z; x_z\le 1\}}. }
\] 
\end{minipage}\hfill\begin{minipage}{0.28\textwidth}
         \includegraphics[page=5,width=\textwidth]{Figures}\end{minipage}\smallskip

When $z+ru_{1:3}\in E'^+$ and
$z\in ([0,1]\times [-r,r])\setminus B^{\cup}(s,t,r)$, we have
$z+ru_{1:3}\in E^+$. This implies that  
\begin{equation*}
0\;\leq\;
\pth{  \indicator{z+ru_{1:3}\in E'^+}%\indicator{z\in [0,1]\times [-r,r]}  
       -  \indicator{z+ru_{1:3}\in E^+} }\indicator{z\in [0,1]\times  [-r,r]}|
\;\leq\;
\indicator{z\in B^{\cup}(s,t,r)}
\end{equation*}
Moreover, we have
\begin{equation*} n^3 \int_{0}^\infty \int_{B^{\cup}(s,t,r)} \int_{\s^3 } \expo{-n\pi r^2}
       \,\cdot\,
         r^3  2\Area{\Delta(u_{1:3})}
        \sigma(du_{1:3}) dz dr 
= 24. 
\end{equation*} This together with \eqref{e:majE} and \eqref{e:majE2} shows that 
\begin{equation*}
{    \frac{64}{3\pi^2} \sqrt{n} -24 \le \;}
\Ex{\sum_{p_{1\diff 3}\in X_n^3} \indicator{\Delta(p_{1:3})\in \DT(X_n)}
       \indicator{p_{1:3}\in E^+}} { \le  \frac{64}{3\pi^2} \sqrt{n} .}
\end{equation*}

Proceeding along the same lines as above, we obtain exactly the same
{bounds} when  we replace $E^+$ by $E^-$. This proves Proposition~\ref{lem:straight} 
by  plugging this value in Equation~(\ref{e:Nst}) since 
\[
\frac{64}{3\pi^2} \sqrt{n} \;-23 
\quad\le\quad \Ex{\frac{N_{s,t} }{\|s - t\|}} 
\quad\le\quad \frac{64}{3\pi^2} \sqrt{n} \;+1. 
\]
\end{prooft}

\section{Length of the Upper Path\label{s:upperBound}}
In this section, we {estimate} the {expectation and the variance}
of the length of the upper path $UP_X$.
{The following proposition deals with the expectation.}

\begin{Prop}\label{lem:upper-bound}
Let $X_n$ be a Poisson point process of intensity $n$ and  $X:=X_n\cup\{s,t\}$ with $s,t\in\R^2$. Let
$\length{UP_{X}}$ be the length of the upper path $UP_X$ in  $\DT(X)$
from $s$ to $t$.
Then 
\[
\Ex{\frac{\length{UP_X}}{\|s - t\|}} =  \frac{35}{3\pi^2} + O\left(n^{-\tfrac{1}{2}} \right) \simeq 1.182.
\]
\end{Prop}
{The above result provides an upper bound for the expectation of the length of the shortest path and implies directly Theorem~\ref{prop:upper-bound}}.

\begin{prooft}
The proof is closely related to the one of
Proposition~\ref{lem:straight}. 
Indeed, assuming without loss of generality that
$s=(0,0)$ and $t=(1,0)$, we define
\begin{equation}\label{def:Lchi}
     L_{X_n}  := \sum_{p_{1\diff 3}\in X_n^3}l_{X_n}(p_{1:3}),
\end{equation}
with 
\[l_{X_n}(p_{1:3}) =  \tfrac{1}{2}\indicator{\Delta(p_{1:3})\in \DT(X_n)}
       \indicator{p_{1:3}\in E^+}
      \| p_2-p_1 \|.\]
In the expression of $ L_{X_n}$, we have considered the lengths of the
edges in $UP_X$, excepted the ones which contain the points $s$ and
$t$. By Proposition~\ref{prop:origin-length}, the expected lengths of
these two edges is $O\pth{n^{-\tfrac{1}{2}}}$. 
Hence, it is
enough to show that 
$\Ex{L_{X_n}} = \frac{35}{3\pi^2} +O\pth{n^{-\tfrac{1}{2}}}$. 
To do it, we apply the Slivnyak-Mecke and
the Blaschke-Petkantschin formulas. This gives
\begin{multline}
\label{eq:Upintegral}
\Ex{L_{X_n}} = \tfrac{1}{2}n^3\!\!\!\int_{\R_+}\!\int_{\R^2}\!\int_{\s^3} \!\!\! \expo{-n\pi r^2}
       \indicator{z+ru_{1:3}\in E^+}
       r \| u_2-u_1 \| \\
         r^3 {{2  \Area{\Delta(u_{1:3})}
      \sigma(du_{1:3})} }dy_z dx_z dr,  
\end{multline}
where we recall that $z=(x_z,y_z)$. 
In the same spirit as Proposition~\ref{lem:straight}, we have:
\begin{multline*}
0\le\;
\pth{  \indicator{z+ru_{1:3}\in E'^+}  -  \indicator{z+ru_{1:3}\in    E^+} }
         \indicator{z\in [0,1]\times  [-r,r]}
         \| u_2 - u_1 \|
\leq 
2\cdot \indicator{z\in B^{\cup}(s,t,r)}
\end{multline*}
and
\begin{align*}
 &n^3 \!\!\! \int_{0}^\infty \!\int_{B^\cup(s,t,r)} \!\int_{\s^3 }\!\!\!
         \expo{-n\pi r^2}
         r^4  { 2 \Area{\Delta(u_{1:3})}
      \sigma(du_{1:3}) }dy_z dx_z dr
\\ & \le 
n^3 \!\!\! \int_{0}^\infty \! 
              \qquad \pi r^2
            \!\int_{\s^3 }\!\!\!
         \expo{-n\pi r^2}
         r^4  { 2 \Area{\Delta(u_{1:3})}
      \sigma(du_{1:3}) }dy_z dx_z dr 
\\ & \le 
\pi n^3\times 
\frac{15}{16\pi^3 n^{\frac{7}{2}}} \times 
24\pi^2= \frac{ 45}{  2\sqrt{n}}
= 22.5\;n^{-\tfrac{1}{2}}. 
  \end{align*}
Replacing $E^+$ by $E'^+$ in (\ref{eq:Upintegral}), it follows that: 

\begin{multline*}
\Ex{L_{X_n}}
  =  
 \tfrac{n^3}{2} \!\!\! \int_{0}^\infty \!\int_{0}^{1}\!\int_{-r}^r \!\int_{\s^3 }\!\!\!
         \expo{-n\pi r^2}
       \indicator{z+ru_{1:3}\in E'^+}
        \| u_2-u_1 \| 
         r^4  \\{{2 \Area{\Delta(u_{1:3})}
      \sigma(du_{1:3})} }dy_z dx_z dr + O\left(n^{-\tfrac{1}{2}} \right).
\end{multline*}
Let $r\geq 0$ be fixed. By taking the change of variables
$h=\frac{y_z}{r}$ (see Figure~\ref{fig:notation}), 
we obtain %for all $x_z\in [-r, 1+r]$
\begin{align*}
\notag &\int_{-r}^r    \int_{\s^3} 
     \indicator{z+ru_{1:3}\in E'^+} \|u_2-u_1\|
      \Area{\Delta(u_{1:3})} \sigma(du_{1:3})dy_z\\  
\notag     & \qquad  =  r \int_{-1}^1 
      \int_{\pi+\smallasin h}^{2\pi-\smallasin h}
     \int_{-\smallasin h}^{\pi+\smallasin h}
      \int_{-\smallasin h}^{\beta_2}
     \det\betamatrix\\
\notag         &     \hspace*{6cm}                \times
              2\sin\frac{\beta_2-\beta_1}{2}
              d\beta_1   d\beta_2d\beta_3 dh
\\ &\qquad = \frac{35\pi}{3}r.
%\label{eq:aireXlength-arcsin}
\end{align*}
It follows that 
\begin{equation*}
\Ex{L_{X_n}}  =  
n^3\cdot \frac{35 \pi}{3}  \int_0^\infty    \int_{0}^{1}      \expo{-n\pi r^2}         r^5 dx_zdr + O\left(n^{-\tfrac{1}{2}} \right)
 % = 
%\tfrac{1}{2}n^3\cdot \frac{70 \pi}{3}  \cdot   \frac{1}{\pi^3 n^3} 
 =  \frac{35}{3\pi^2}   + O\left(n^{-\tfrac{1}{2}} \right).
\end{equation*}
\end{prooft}

{The following proposition deals with the variance of the length of the upper path.}

\begin{Prop}\label{prop:up-variance}
Let $X_n$ be a Poisson point process of intensity $n$ and  $X:=X_n\cup\{s,t\}$ with $s,t\in\R^2$. Let
$\length{UP_{X}}$ be the length of the upper path $UP_X$ in  $\DT(X)$
form $s$ to $t$. Then 
\[\VVV{\frac{\ell(UP_{X})}{\|s-t\|}} = O\pth{n^{-\tfrac{1}{2}}}.\]
\end{Prop}

The proof uses similar tools as the one for the expected value and is
postponed to Appendix~\ref{ap:variance}. 
{As a corollary, we obtain an estimate of the tail of the length of the upper path.}

\begin{Cor}\label{cor:sp-is-small}
{With the same notation as above, we have}
\[
\Prob{\ell(UP_{X}) > 1.2 ||s-t|| } 
\quad = \quad 
O\pth{n^{-\tfrac{1}{2}}}.
\]
\end{Cor}
\begin{prooft}
It follows from the Chebyshev's inequality that
\begin{equation*}
\Prob{\ell(UP_{X}) > 1.2 ||s-t|| }  \leq \frac{\VVV{\ell(UP_{X})}}{1.2\|s-t\| - \Ex{    \ell(UP_{X}) }}.
 \end{equation*} 
 {This concludes the proof according to}  Proposition~\ref{lem:upper-bound}  and   Proposition~\ref{prop:up-variance}.
\end{prooft}

\section{Lower Bound on Shortest Path \label{s:lowerBound}} 

In this section, we prove Theorem~\ref{th:lower}. By scaling
  invariance, our problem is the same as if the intensity $n$ is
  constant and $s=(0,0)$, $t=(k,0)$, where $k\in\N^*$ goes to
  infinity.  

Before starting our proof, we give the main ideas. First, we
discretize the plane $\R^2$ into squares $C(v)$, which will be called
pixels, with $v\in \Z^2$. The set of pixels intersecting the smallest
path $SP_X$ will be called the lattice animal associated with
$SP_X$. In a sense which will be specified, we associate with each
pixel $C(v)$ a (strong) horizontality property which ensures that there exists a
path "almost horizontal" in the Delaunay triangulation through the
pixel $C(v)$. Then, we proceed into four steps.  
\begin{enumerate}
\item We begin with preliminaries by introducing formally the notions
  of animal lattice and strong horizontality property.
\item We give an upper bound for the size of the associated
  lattice animal (Corollary~\ref{cor:size-animal-proba}). 
Moreover, given a property $\mathcal{Y}$, we provide
  an upper bound for the probability that the number of pixels in the
  lattice animal with property $\mathcal{Y}$ is large. This result
  will be applied to the case where $\mathcal{Y}$ is the horizontality
  property (Proposition~\ref{lem:percolation}). 
\item We establish a lower bound for the length of any path with respect
  to the number of pixels with a strong horizontality property in the associated lattice animal (Lemma~\ref{lem:length-animal}). 
\item We prove that the probability that a pixel has a (strong) horizontality property is small (Proposition~\ref{lem:good-proba}).
\item We conclude our proof by combining the steps described above (Section~\ref{sec:wrapup}).
\end{enumerate}

\subsection{Preliminaries}
\subsubsection{Animal lattices}
\begin{figure}[t]
     \begin{center}
         \includegraphics[page=7,width=\textwidth]{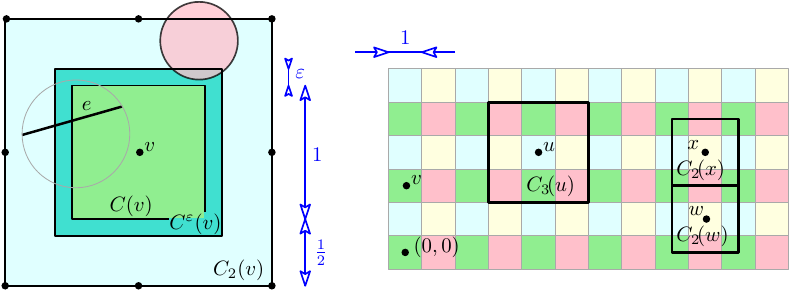}
     \end{center}
     \caption{Squares definitions.
         %   {\bf Left:}  The squares $C(v)$ and $C_2(v)$ with auxiliary
         % constructions for the proofs of Lemmas~\ref{lem:choose-grid}
         % and~\ref{lem:horizontality}.
          %{{\bf Right:}          Coloring of pixels.}
        \label{fig:pixel}
     }
 \end{figure}

Recall that $s=(0,0)$ and $t=(k,0)$ for some integer $k$.
We discretize $\R^2$ into pixels as follows. 
Let $\grid=(\Z^2,E)$ be the graph with set of edges satisfying $(v,w)\in E
\Leftrightarrow ||v-w||=1$, where $v,w\in\Z^2$.  In
digital geometry, the graph is known as a {\em 4-connected neighborhood} and
each vertex of $\grid$ is called a {\em pixel}. 
Moreover, for each $v\in\Z^2$, we
consider different scaled versions  of squares (see: Figure~\ref{fig:pixel}-left)
\[C(v):=v\oplus[-\tfrac{1}{2},\tfrac{1}{2}]^2, \quad
C^\varepsilon(v):=v\oplus[-\tfrac{1}{2}-\varepsilon,\tfrac{1}{2}+\varepsilon]^2;
\quad C_\lambda(v):=v\oplus[-\tfrac{\lambda}{2},\tfrac{\lambda}{2}]^2,\]
where $\oplus$ denotes the Minkowski
sum, $\lambda\in\Z$ and $\varepsilon\in\R^*_+$. With a slight abuse of notation, we also say
that $C(v)$ is a pixel. 

We define scaled and translated
versions of the grid $\grid$ as follows. For $\lambda\in\Z$ and $\tau\in\Z^2$
we denote by $\lambda \grid$ the grid of points in $(\lambda \Z)^2$
with edges of length $\lambda$
and $\tau+\lambda \grid$ its translation by $\tau$.
We also split $\grid$ in 4 subgrids as follows: each subgrid is referred to {as} a color $c\in Colors$, where $Colors:=\{green, pink, blue, yellow\}$. Each subgrid with color $c$ is denoted:
\[ \grid_c=2\grid+O_c,\]
where 
\[ 
O_{(green)}=(0,0),\quad 
O_{(pink)}=(1,0),\quad 
O_{(blue)}=(0,1), \quad \mbox{and } 
O_{(yellow)}=(1,1).
\]

Notice that when $v$ and $w$ are two pixels with the same color, the squares $C_2(v)$ and $C_2(w)$ have disjoint interior. Such a property will be useful to ensure the independence of a suitable family of random variables.
 
\begin{figure}[t]
     \begin{center}
         \includegraphics[page=6,scale=1]{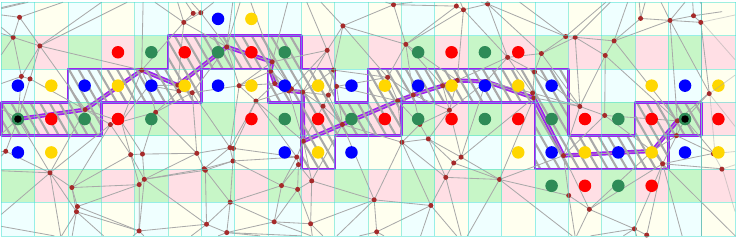}
     \end{center}
     \caption{
       Lattice animals. The pixels $C(v)$ for   $v\in\animal(P)$ are hashed
       and the path $P$ is purple.
       The  $\animal_{(c)}(P)$, for each $c\in Colors$,
       are {referred to as} big dots of the relevant color.
        \label{fig:grid-colors}
     }
 \end{figure}

We conclude this section with the so-called notion of animal~\cite{cox93}. Given a graph $G':=(V',E')$, a {\em lattice animal}
is a collection of vertices $\animal\subset V'$ such that for every pair of distinct vertices $v,w\in
\animal$ there is a path in {$G'$} connecting $v,w$ visiting only
vertices in $\animal$.  With each path $P$ in the Delaunay triangulation, we associate the so-called lattice animal of $P$ in $\grid$ (see Figure~\ref{fig:grid-colors}):
\[\animal(P) = \{ v \mbox{ vertex of }\grid: \,C(v)\cap P\neq\varnothing\}.\] 
 In the same spirit as above, for each $c\in Colors$, we let 
\[\animal_{(c)}(P) = \{ v \mbox{ vertex of }\grid_{(c)}: \,C_2(v)\cap
P\neq\varnothing\}\] 
We also considered animals at different scales (see Figure~\ref{fig:percolation}):
\[\animal_{(c),\lambda}(P) = \{ v \mbox{  vertex of }\lambda\grid+O_{c}: \,C_\lambda(v)\cap P\neq\varnothing\}.\] 

\subsubsection{Properties on pixels}
For any pixel $v\in\Z^2$, we 
consider three events, namely ${\cal I}_\varepsilon(v) $, ${\cal H}_\rho (v) $,  and 
${\cal  H}'_{\varepsilon,\alpha,\kappa}(v) $. These events depend on four parameters $\varepsilon>0$,
$\rho>0$, $\alpha\in [0,\frac{\pi}{2}]$, and $\kappa>1$ and  are described 
below. 
\begin{itemize}
\item \textit{Independence property} ${\cal I}_\varepsilon(v)$: this event holds if, for  
any Delaunay triangle in $\DT(X)$ intersecting $C^\varepsilon(v)$, the circumdisk {of the triangle is included in} $C_2(v)$ and  $\|s-v\|\ge 2$ 
and $\|t-v\|\ge 2$ (to avoid pixels that are neighbors of $s$ or $t$). {Notice that this event is $\sigma(X_n\cap C_2(v))$ measurable}.
\item \textit{Strong horizontality property}  ${\cal H}_\rho(v)$: this event holds if there exists a path along Delaunay edges in $\DT(X)$, between $x_v-\tfrac{1}{2}$ and $x_v+\tfrac{1}{2}$, {intersecting} 
$C(v)$ and with length smaller than $1+\rho$. A path satisfying this
property is denoted by ${\cal PH}_\rho(v)$. Besides, the first
and the last edges of such a path are clipped by the vertical 
lines $x=x_v-\tfrac{1}{2}$ and $x=x_v+\tfrac{1}{2}$.

\item \textit{Weak horizontality property} ${\cal
    H}'_{\varepsilon,\alpha,\kappa}(v)$: this event holds if the total
  length $L_{\varepsilon,\alpha,\kappa}(v) $ of the horizontal
  projection of all edges in $\DT(X)$ intersecting $C^\varepsilon(v)$
  and having an angle with respect to
the $x$-axis smaller than $\alpha$, is greater than $\tfrac{1}{\kappa}$, i.e. 
\[
L_{\varepsilon,\alpha,\kappa}(v) = \sum_{e, e\cap C^\varepsilon(v)\neq \varnothing, \eangle{e}\le \alpha} h(e) \ge \tfrac{1}{\kappa}.
\]
\end{itemize}

\subsection{Animals\label{sec:animals}}

\subsubsection{Size of animals\label{sec:animal-size}}

We establish below a series of results on the size of animals
obtained from a discretization of some path starting 
from $s$ and going to $t$ with $s,t\in\Z^2$.
The following lemma is due to  Gerard, Favreau and
Vacavant~\cite{gerard2015tight}. For completeness, 
we give a more concise proof of their result.  

\begin{Le}\label{lem:size-animal}
Let $s,t\in\R^2$ and let $P$ be a path between $s$ and $t$ (not necessarily in $\DT{X}$). Then
\[\size{\animal(P)}\leq \tfrac{3\sqrt{2}}{2}\length{P}+1\simeq 2.12 \length{P}+1. \]
\end{Le}

\begin{figure}[t]
     \begin{center}
           \includegraphics[page=8,width=\textwidth]{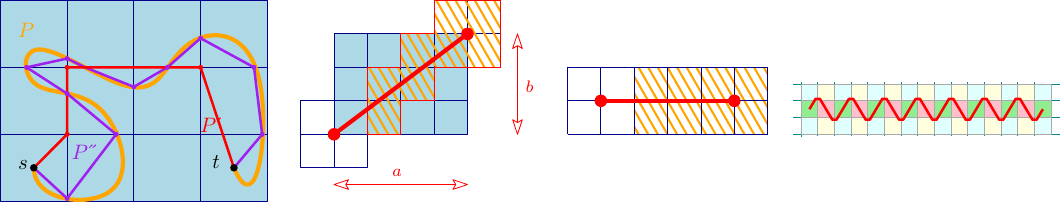}
     \end{center}
     \caption{
Illustration for the proof of Lemma~\ref{lem:size-animal}. 
Left: A realization of the paths $P$, $P''$, and $P'$.
Center: hatched pixels are
charged to the red edge.
Right: The red path visits 3 pixels per column and has a length
       close to $\sqrt{2}$ per column (2 pixels in the first and last
       column with length $\frac{\sqrt{2}}{2}$ in these columns).
        \label{fig:numbpixels}
     }
 \end{figure}

\begin{prooft}
Let $P:=\{e_1, \ldots, e_k\}$ be a path between $s$ and $t$, with
$e_i:=[Z_{i-1}, Z_{i}]$ for each $1\leq i\leq k$ and $Z_0=s$ and
$Z_k=t$. 
We can construct an auxiliary path $P''$ with the same lattice animal
as $P$ and such that $\ell(P'')\leq \ell(P)$ by straightening the path
$P$ between 
the intersections of $P$ and the boundaries of pixels $C(v)$, $v\in\Z^2$
(see Figure \ref{fig:numbpixels}-left).
 Then, we 
can construct a path $P'$ from $P''$, again with the same lattice
animal as $P$, with  $\ell(P')\leq \ell(P'')$, and such that the
vertices of $P'$ belong to 
the set of corners of pixels. The path $P'$ is defined by moving the vertices of $P''$ along
 the boundaries of pixels in the direction that shorten the path, up to a corner or up to 
the point that aligns the two incident segments. In particular, the new path is such that
$\animal(P)=\animal(P')$ and $\ell(P)\geq \ell(P')$. Hence, without
loss of generality, we can assume that 
$Z_i\in(\tfrac{1}{2},\tfrac{1}{2})+\Z^2$ for each $1\leq i\leq k-1$. 

Denoting by $E_i:=\bigcup_{j=0}^i e_j$, we trivially obtain that
\[\size {\animal(P)} = \sum_{i=1}^k \size{ \animal'(e_i)},\]
where 
\[\animal'(e_i)=\{v\in \Z^2: C(v)\cap e_i\neq \varnothing \text{ and } C(v)\cap E_{i-1} = \varnothing\}.\] 
We provide below an upper bound for $\size{ \animal'(e_i)}$.
Indeed, for each edge $e_i:=[Z_{i-1}, Z_i]$,
$2\leq i\leq k-1$, 
we consider the rectangle 
$[x_{Z_{i-1}}, x_{Z_i}]\times [y_{Z_{i-1}}, y_{Z_i}]$
(see Figure \ref{fig:numbpixels}-center).
If $x_{Z_{i-1}}=x_{Z_i}$ or $y_{Z_{i-1}}=y_{Z_i}$, we trivially obtain that $\size{ \animal'(e_i)}\leq 2\cdot \ell(e_i)$. If not, we notice that $\size{ \animal'(e_i)}=a_i+b_i+1$ and $\ell(e_i)=\sqrt{a_i^2+b_i^2}$ with $a_i=|x_{Z_{i-1}} - x_{Z_i}|$ and $b_i=|y_{Z_{i-1}} - y_{Z_i}|$. Besides, we can easily prove that for each $a,b>0$, we have $\tfrac{a+b+1}{\sqrt{a^2+b^2}}\leq \tfrac{3\sqrt{2}}{2}$ 
(the equality holds when $a=b=1$). This implies that $\size{
  \animal'(e_i)}\leq \tfrac{3\sqrt{2}}{2}\cdot \ell(e_i)$. 
Finally we have $\ell(e_1)=\ell(e_k)=\tfrac{\sqrt{2}}{2}$, $\size{\animal'(e_1)}=4$, and $\size{\animal'(e_k)}=0$.
It follows that
\[\size {\animal(P)} \leq  4+ \sum_{i=2}^{k-1} \tfrac{3\sqrt{2}}{2}\cdot \ell(e_i) = \tfrac{3\sqrt{2}}{2}\cdot \ell(P)+1.\]
In Figure \ref{fig:numbpixels}-right, we depict a path $P$ such that 
$\size{\animal(P)}\simeq \tfrac{3\sqrt{2}}{2}\length{P}+1\simeq 2.12\length{P}+1$, 
which shows that our upper bound in Lemma~\ref{lem:size-animal} is tight.

\end{prooft}

As a corollary, we obtain a deterministic upper bound for the size of the animal associated with the shortest path.

\begin{Cor}\label{cor:size-animal}
Let $SP_X\in\setpath$ with $\|s-t\|=k$. Then 
\[\size{\animal(SP_X)}< 4.24\,k+1 \qquad\mbox{ and}\qquad \forall \lambda\in\Z,\;\size{\animal_{(c),\lambda}(SP_X)}< \frac{4.24\,k}{\lambda}+1.\]
\end{Cor}
\begin{prooft}
The first result is a direct consequence of Lemma~\ref{lem:size-animal} and the fact that 
$\length{SP_X}<1.998 k$ according to the main result of Xia~\cite{x-sfdtl-13}.
The second one is obtained by applying the first result to the grid $\lambda\grid$.
\end{prooft}
The following result provides a better upper bound for the size of the animal with high probability.
\begin{Cor}\label{cor:size-animal-proba}
Let $SP_X\in\setpath$ with $\|s-t\|=k$. For any $\lambda>0$, $c\in Colors$, and $k>0$, let $\mathcal{E}(c,k, \lambda)$ be the event:
\begin{equation}
\label{eq:defeventE}
\mathcal{E}(c,k, \lambda):=\left\{ \size{\animal_{(c),\lambda}(SP_X)}< \tfrac{2.55\,k}{\lambda}+1 \right\}. \end{equation}
 Then $\Prob{\mathcal{E}(c,k, \lambda)} =  1-O\left( k^{-\tfrac{1}{2}}\right)$.
\end{Cor}
\begin{prooft}
This is a probabilistic version of Corollary~\ref{cor:size-animal}
obtained by bounding the length of $SP_X$ by the one of $UP_X$
instead of using the deterministic bound due to Xia.
Corollary~\ref{cor:sp-is-small}
ensures that the upper bound for the length of $UP_X$ holds with
probability $1-O\left(k^{-\tfrac{1}{2}}\right)$.
The constant comes from the fact that $1.2\cdot 2.12 < 2.55$.
\end{prooft}

\subsubsection{Animal lemma}
In this section, we establish a result which ensures that the number
of pixels in $\animal(SP_X)$ which satisfy a property $\mathcal{Y}$ is
not large with high probability. In the sequel, we will use the
following notation: 
\[\Z^2_{s,t}:= \Z^2\setminus\{v\in \Z^2: \|v-s\|\leq 2 \text{ or } \|v-t\|\leq 2\}.\]
For each color $c$, we also let $\Z^{2, (c)}_{s,t}:=\Z^2_{s,t}\cap \{v\in \Z^2 \text{ with color c}\}$.

\begin{Prop}
\label{lem:percolation}
Let $p\in(0,0.01]$ and let $\mathcal{Y}:=(Y_v)_{v\in\Z^2}$ be a family
of events such that, for any color $c$, the events 
$(Y_v)_{v\in\Z^{2,(c)}_{s,t}}$ are independent and $p=\Prob{Y_v}$ for
each $v\in\Z^{2,(c)}_{s,t}$. 
For any $\mathbf{A}\subset\Z^2$, we denote by  
$\sharp_{\cal  Y}(\animal)=\sum_{v\in \animal} \indicator{Y_v}$  the number of
pixels $v$ in $\animal$ such that the event $ Y_v$ holds. Then, we have 
    \begin{align*}
       \Prob{ \sharp_{\cal Y} (\animal_{(c)}(SP_X))  \ge 4 k\sqrt{p}} =    O\pth{k^{-\frac{1}{2}}}.
    \end{align*}
\end{Prop}

\begin{prooft}
Our proof relies on an adaptation of a result due to  Devillers and
Hemsley~\cite[Lemma~7]{devillers:hal-01216212}. The main idea is to
discretize the shortest path at different scales and to use standard
ideas of (site) percolation theory.  
Let $\lambda\in 4\Z+2$. This restriction on $\lambda$ ensures that 
$\forall v\in \grid_c,\;
\forall w\in\lambda\grid+O_c :\;
C_2(v)\cap C_\lambda(w)\neq\varnothing \Rightarrow C_2(v)\subset
C_\lambda(w)$. 
Taking $\mathcal{E}(c,k, \lambda)$ as in \eqref{eq:defeventE}, we obtain for each $x>0$, that 
\begin{multline*}
    \Prob{ \sharp_{\cal Y} (\animal_{(c)}(SP_X))  \ge  x k \sqrt{p}} = \Prob{ \{\sharp_{\cal Y} (\animal_{(c)}(SP_X))  \ge  x k \sqrt{p} \} \cap  \mathcal{E}(c,k, \lambda)} \\
  +   \Prob{ \{\sharp_{\cal Y} (\animal_{(c)}(SP_X))  \ge  x k \sqrt{p}\}  \cap   (\neg  \mathcal{E}(c,\lambda,k))}.
  \end{multline*}
From Corollary~\ref{cor:size-animal-proba}, is follows that
\begin{equation*}
 \Prob{ \sharp_{\cal Y} (\animal_{(c)}(SP_X))  \ge  x k \sqrt{p}} \leq    \Prob{ \{\sharp_{\cal Y} (\animal_{(c)}(SP_X))  \ge  x  k \sqrt{p}\}\cap  \mathcal{E}(c,k,\lambda) }
         + O\pth{k^{-\tfrac{1}{2}}}.
\end{equation*}
 The  
animal $\animal_{(c),\lambda}(SP_X)$ can be viewed as a sequence of connected squares of size
$\lambda$, starting at the square containing $s$. Besides, on the
event $\mathcal{E}(c,k,\lambda)$, the animal
$\animal_{(c),\lambda}(SP_X)$ belongs to the family ${\cal
  A}_{(c),\lambda}(k)$ of animals $\animal\subset\lambda\grid +O_c$
such that $O_c\in\animal$ and
$\size{\animal}\leq\floor{\frac{2.55k}{\lambda}+1}$. 
Each animal $\mathbf{A}\in {\cal A}_{(c),\lambda}(k)$ can be encoded
as a word on four letters $\{S,E,N,W\}$ (standing for south, east,
north, and west) of length $\size{\animal}$. Hence, 
$\size{ {\cal  A}_{(c),\lambda}(k) }
\leq
4^{\floor{\frac{2.55\,k}{\lambda}+1}}
\leq
4^{{\frac{2.55\,k}{\lambda}+1}}$. 
Moreover, since $\lambda\in 4\Z+2$, we have   
\[
\animal_{(c)}(SP_X) 
\subset \animal_{(c),\lambda}^{(\lambda)}(SP_X),
\]
where, for all animal $\mathbf{A}\subset \lambda\grid+O_c$, we let $\mathbf{A}^{(\lambda)}:=\{v\in
\grid_{(c)}: \exists w\in \mathbf{A}, v\in C_\lambda(w)\}$
(see Figure~\ref{fig:percolation}).
\begin{figure}[t]
     \begin{center}
           \includegraphics[page=9,width=0.7\textwidth]{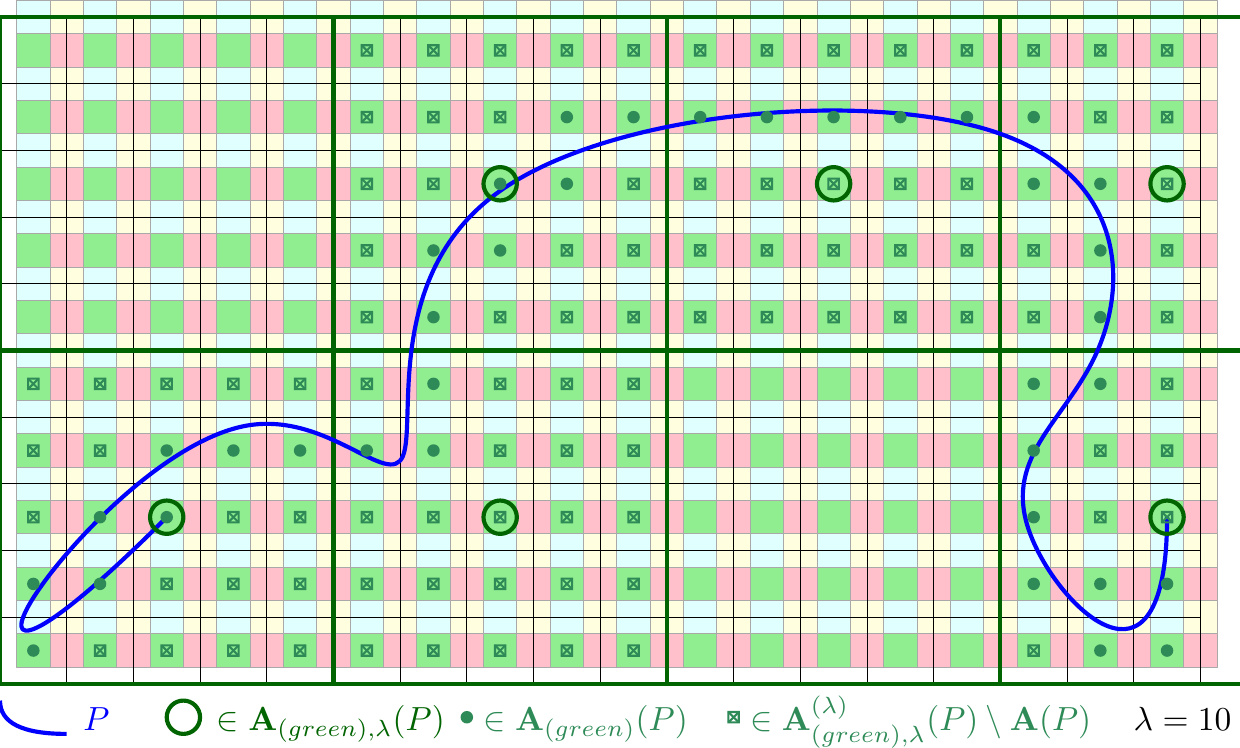}
     \end{center}
     \caption{The set $\animal_{(green),\lambda}^{(\lambda)}(P)$
         for a path $P$ .
        \label{fig:percolation}
     }
 \end{figure}
This implies that for each $x>0$, we have 
\begin{multline}\label{eq:boundbinomial}
\Prob{ \{\sharp_{\cal Y} (\animal_{(c)}(SP_X))  \ge  x k \sqrt{p}\}\cap \mathcal{E}(c,k, \lambda)}\\
\begin{split}
& \leq \Prob{ \{\sharp_{\cal Y} (\animal_{(c),\lambda}^{(\lambda)}(SP_X))  \ge  x k \sqrt{p}\}
           \cap  \mathcal{E}(c,k, \lambda)}\\ 
& \leq \Prob{ \bigcup_{\animal\in \mathcal{A}_{(c),\lambda}(k)} \left\{ \sharp_{\cal Y}  (\animal^{(\lambda)})  \ge x  k \sqrt{p}\right\} }\\ 
& \leq \sum_{\animal\in \mathcal{A}_{(c),\lambda}(k)}
         \Prob{ \sharp_{\cal Y} (\animal^{(\lambda)})  \ge x  k \sqrt{p} } \\
& \leq\sum_{\animal\in \mathcal{A}_{(c),\lambda}(k)}
       \Prob{\sum_{v\in \animal^{(\lambda)} \cap \Z^{2,(c)}_{s,t}}\indicator{Y_v}\geq x k \sqrt{p}-10   }.
\end{split}  
  \end{multline}
The last inequality comes from the fact that 
$\sum_{v\in  \animal^{(\lambda)} \cap \Z^{2,(c)}_{s,t}}\indicator{Y_v}\geq \sharp_{\cal Y} (\animal^{(\lambda)}) -10$ 
since  the number of pixels in $\Z^2$ with color $c$, which are not in $\Z^{2,(c)}_{s,t}$, is lower than 10.
From the assumption of Proposition~\ref{lem:percolation}, 
we know that for each $\mathbf{A}\in \mathcal{A}_{(c),\lambda}(k)$,
the random variable  
$\sum_{v\in  \animal^{(\lambda)} \cap \Z^{2,(c)}_{s,t}}\indicator{Y_v}$ 
is a binomial distribution with parameters $(N,p)$, where 
$N := \size{\animal^{(\lambda)}\cap  \Z^{2,(c)}_{s,t}}\leq  0.64\lambda k$
since
$\size{\animal^{(\lambda)}}\le\frac{2.55k}{\lambda}\frac{\lambda^2}{4}+1$,
$s\in \animal^{(\lambda)}\setminus  \Z^{2,(c)}_{s,t}$, and
$\frac{2.55}{4}\leq 0.64$.

According to the Chernoff's
inequality, we obtain for all $\animal\in\mathcal{A}_{(c),\lambda}(k)$ that 
\begin{eqnarray*}
\Prob{\sum_{v\in \animal^{(\lambda)} \cap    \Z^{2,(c)}_{s,t}}\indicator{Y_v}\geq x k \sqrt{p}-10   } 
& \leq &
e^{-x k \sqrt{p}+10} \cdot \Ex{  e^{\sum_{v \in \animal^{(\lambda)} \cap  \Z^{2,(c)}_{s,t}}   \indicator{Y_v}}}
\\ & = &
e^{-x k \sqrt{p}+10} \cdot  ( (1-p) + p e)^{\size{\animal^{(\lambda)}\cap  \Z^{2,(c)}_{s,t} }}
\\ & \leq &
e^{-x k \sqrt{p}+10} \cdot  ( 1+p(e-1))^{0.64\lambda k}.
\end{eqnarray*}
This together with \eqref{eq:boundbinomial} and the fact that $\size{\mathcal{A}_{(c),\lambda}(k)} \leq 4^{{\frac{2.55k}{\lambda}+1}}$ implies that 
\begin{equation*}
\Prob{ \{\sharp_{\cal Y} (\animal_{(c)}(SP_X))  \ge  x  k \sqrt{p}\} \cap   \mathcal{E}(c,k,\lambda) } 
\leq 4^{{\frac{2.55k}{\lambda}+1}} \cdot 
e^{-x k \sqrt{p}+10} \cdot  ( 1+p(e-1))^{0.64\lambda k}. 
\end{equation*}

With standard computations, we get 
\begin{multline*}
\Prob{ \{\sharp_{\cal Y} (\animal_{(c)}(SP_X))  \ge  x  k \sqrt{p}\} \cap   \mathcal{E}(c,k,\lambda) }
\\
\begin{split}
& \leq 4^{{\frac{2.55k}{\lambda}+1}} \cdot 
e^{-x k \sqrt{p}+10} \cdot  ( 1+p(e-1))^{0.64\lambda k} 
\\ &  \leq 
4^{{\frac{2.55k}{\lambda}+1}} \cdot e^{-x k \sqrt{p}+10} \cdot  ( e^{p(e-1)})^{0.64\lambda k} 
\\ & = 
4e^{10}\cdot \exp\pth{\log 4\cdot {\frac{2.55k}{\lambda}} -x  k \sqrt{p} +0.64\ p(e-1)\lambda k }
\\ & \leq 
4e^{10} \cdot \exp\pth{-k\sqrt{p}\pth{x -\pth { \frac{3.54}{\lambda\sqrt{p}} + 1.1 \lambda\sqrt{p}}}}
\\  & \leq 
4e^{10}\ e^{-(x-3.98) k\sqrt{p}}\qquad\mbox{by choosing   $\lambda\in\left[\frac{1.6}{\sqrt{p}},\frac{2}{\sqrt{p}}\right]\cap(4\Z+2)$}
\\ & = 
o\pth{k^{-\frac{1}{2}}}\qquad\mbox{by choosing $x=4$ for any fixed $p$ when $k\rightarrow\infty$}.
\end{split}
\end{multline*}
The fact that a suitable
$\lambda\in\left[\frac{1.6}{\sqrt{p}},\frac{2}{\sqrt{p}}\right]\cap(4\Z+2)$
exists is ensured by the condition $p<0.01$.
\end{prooft}

\subsection{A lower bound for the length of a path}
We establish below a lower bound for the length of any path with respect to the number of pixels with a strong horizontality property.

\begin{Le}\label{lem:length-animal}
Let $P\in\setpath$ and $\rho>0$. Then 
\[
\length{P} \ge k + \rho  \pth{k - 4 \max_{c\in Colors}\sharp_{\cal H}(\animal_{(c)}(P)) },
\]
where $\sharp_{\cal H}(\animal_{(c)}(P))$ is the number of  pixels 
$v$ in $\animal_{(c)}(P)$ such that ${\cal H}_\rho(v) \vee \neg {\cal I}_{\varepsilon_\rho}(v)$.
\end{Le}
\begin{prooft}
Splitting the path $P$ into vertical columns, we have $\length{P} = \sum_{i\in\Z}\length{P\cap Col[i]},$ 
where, for each $i\in\Z$,
$Col[i]=[i-\tfrac{1}{2},i+\tfrac{1}{2}]\times\R$ is the $i^th$ column.
For $i\in\{1,\ldots, k-1\}$,  let $v[i]$ be the lowest pixel of  $\animal(P)\cap Col[i]$
such that there is a connected component of $P\cap Col[i]$
intersecting $C(v[i])$ and the left and right side of $Col[i]$.
Notice that such a pixel exists for each column {$1\le i\le k-1$.}
On the event $\neg{\cal H}_\rho(v[i])  \wedge {\cal  I}_{\varepsilon_\rho}(v[i])$, we use the fact that 
$\length{P\cap Col[i]}\geq 1+\rho$. 
On the complement of this event,  we use the trivial inequality $\length{P\cap Col[i]}\geq 1$.  
Denoting by $N:={\sum_{i=1}^{k-1}}\indicator{ {\cal H}_\rho(v[i])  \vee \neg {\cal I}_{\varepsilon_\rho}(v[i])}$ the number of horizontal pixels on the path, it follows that
\begin{align*}
\length{P} 
& \geq  
{\sum_{i=1}^{k-1}}\indicator{ {\cal H}_\rho(v[i])  \vee \neg {\cal I}_{\varepsilon_\rho}(v[i])}
+{\sum_{i=1}^{k-1}} (1+\rho)\indicator{\neg{\cal H}_\rho(v[i])  \wedge {\cal  I}_{\varepsilon_\rho}(v[i])}
\\& \geq N + (1+\rho)(k-N) = k + \rho (k-N).
 \end{align*}
Then we conclude the proof by observing that
\[
N
%=\sum_{i=0}^k\indicator{ {\cal H}_\rho(v[i])  \vee \neg {\cal I}_{\varepsilon_\rho}(v[i])}
 \leq  \!\!\sum_{v\in\animal(P)}\!\!\indicator{ {\cal H}_\rho(v)  \vee \neg {\cal I}_{\varepsilon_\rho}(v)}
\leq \!\!\!\! \sum_{c\in Colors} \!\!\!\!\sharp_{\cal H}(\animal_{(c)}(P))
\leq  \!\!4\!\!\!\! \max_{c\in Colors} \!\!\!\!\sharp_{\cal H}(\animal_{(c)}(P)). 
\]
\end{prooft}

\subsection{Probability that a pixel has a strong horizontality property}
To apply Proposition~\ref{lem:percolation}, we have to estimate the probability of the event $  {\cal H}_{\rho}(v) \vee \neg{\cal I}_{\varepsilon_\rho}(v) $. An upper bound for this probability is given in the following result.

\begin{Prop} 
\label{lem:good-proba}
Let $v\in\Z^2_{s,t}$   and $0<\rho<4\cdot 10^{-6}$  
and let $\varepsilon_\rho:=\sqrt{\rho}\sqrt{2+\rho}$.  
 Then
\[ \Prob{    {\cal H}_{\rho}(v) \vee \neg{\cal I}_{\varepsilon_\rho}(v)}
\le   
P(\rho,n),
\]
where 
\[P(\rho,n):= 95 n^3 \expo{-0.194  n} + \left(19n^2+13 n +4  \right)e^{-n\pi} % lemma independence
+ 
31.76\pth{\frac{3}{4}+\frac{\sqrt{\rho}\sqrt{\rho+2}}{2}}^2 \sqrt{\rho n}  % lemma horizontal
.\]
\end{Prop}
One of the main difficulties to prove Theorem~\ref{th:lower} is actually contained in the above result. To pave the way, we proceed into two steps. First, we choose  parameters $\varepsilon$ and $\alpha$ in such a way that the strong horizontality property is stronger than the weak horizontality property. Secondly, we provide bounds for the probability that a pixel has a weak horizontality property or an independence property.

\subsubsection{Strong vs weak horizontality}

\begin{Le}\label{lem:epsilon}
Let $v\in\Z^2$ and $\rho>0$ 
and let $\varepsilon_\rho:=\sqrt{\rho}\sqrt{2+\rho}$.  
 If the property ${\cal H}_\rho(v)$ holds, then  ${\cal PH}_\rho(v)\subset C^{\varepsilon_\rho}(v)$. 
\end{Le}
\begin{prooft}
Assume that $v=0$ without loss of generality. Up to a vertical
translation, the shortest path between the lines $x=-\tfrac{1}{2}$ and
$x=\tfrac{1}{2}$ crossing $C(0)$ and intersecting the {complement} of
$C^{\varepsilon_\rho}(0)$ is the segment from
$\pth{-\tfrac{1}{2},\tfrac{1}{2}}$ to
$\pth{\tfrac{1}{2},\tfrac{1}{2}+\varepsilon_\rho}$. Since the length
of this segment is $\sqrt{1+\varepsilon_\rho^2}=1+\rho$, we
necessarily have ${\cal PH}_\rho(v)\subset C^{\varepsilon_\rho}(v)$.  
\end{prooft}

\begin{Le}\label{lem:alpha}
Let $v\in\Z^2$ and $\rho>0$ and $0<\kappa<1$ such that
$\frac{\kappa}{\kappa-1}\rho<\frac{\pi^2}{8}$ and let 
$\alpha_{\rho,\kappa}:=\sqrt{2\frac{\kappa}{\kappa-1}\rho}$. If the property ${\cal H}_\rho(v)$ holds
then the same is true for 
${\cal  H}'_{\varepsilon_\rho,\alpha_{\rho,\kappa},\kappa}(v)$. 
 \end{Le}
\begin{prooft}
We make a proof by contradiction, assuming
${\cal H}_\rho(v)$ and $\neg{\cal  H}'_{\varepsilon_\rho,\alpha_{\rho,\kappa},\kappa}(v)$. 
The main idea is to split
the edges $e$ in  ${\cal PH}_\rho(v)$ with respect to their angles
$\eangle{e}$. Indeed,  
\[1+\rho \ge \length{ {\cal PH}_\rho(v) }  = \sum_{e\in {\cal
    PH}_\rho(v), \eangle{e}\le \alpha_{\rho,\kappa}} \length{e} + \sum_{e\in
  {\cal PH}_\rho(v), \eangle{e}> \alpha_{\rho,\kappa}} \length{e},\] where the
inequality comes from the property ${\cal H}_\rho(v)$. For each $e\in
{\cal PH}_\rho(v)$, we use the trivial inequality $l(e)\geq h(e)$ when
$\eangle{e}\le \alpha_{\rho,\kappa}$. If $\eangle{e}> \alpha_{\rho,\kappa}$, we notice
that,
for $\frac{\kappa}{\kappa-1}\rho<\frac{\pi^2}{8}$,
\[
l(e)
>     \frac{h(e)}{\cos\alpha_{\rho,\kappa}}
\ge \pth{1+\frac{\alpha_{\rho,\kappa}^2}{2}}h(e) 
= \pth{1+\frac{\kappa}{\kappa-1}\rho}h(e),
\] where the second inequality comes from the fact that $\frac{1}{\cos\alpha}
\ge \pth{1+\frac{\alpha^2}{2}}$ for any $\alpha \in [0,\tfrac{\pi}{2})$.

 It follows that 
\begin{align*} 
1+\rho & 
> \pth{1+\frac{\kappa}{\kappa-1}\rho-\frac{\kappa}{\kappa-1}\rho} \sum_{e\in {\cal PH}_\rho(v), \eangle{e}\le \alpha_{\rho,\kappa}}\!\!\!\!h(e) 
+ \pth{1+\frac{\kappa}{\kappa-1}\rho}  \sum_{e\in {\cal PH}_\rho(v), \eangle{e}> \alpha_{\rho,\kappa}}\!\!\!\!h(e)\\
& = \pth{1+\frac{\kappa}{\kappa-1}\rho} \sum_{e\in {\cal PH}_\rho(v)}h(e)\quad-\frac{\kappa}{\kappa-1}\rho\sum_{e\in {\cal PH}_\rho(v), \eangle{e}\le \alpha_{\rho,\kappa}} h(e).
\end{align*}
By assumption, the property 
${\cal H}'_{\varepsilon_\rho,\alpha_{\rho,\kappa},\kappa}(v)$  
does not hold. 
Then we deduce from Lemma~\ref{lem:epsilon} that 
\[
1+\rho 
> \pth{1+\frac{\kappa}{\kappa-1}\rho} - \frac{\kappa}{\kappa-1}\rho L_{\varepsilon_\rho,\alpha_{\rho,\kappa},\kappa}(v)
> 1+\frac{\kappa}{\kappa-1}\rho - \frac{\kappa}{\kappa-1}\rho \frac{1}{\kappa}
=1+\rho
,\]
getting a contradiction.
\end{prooft}

\subsubsection{Pixel probabilities}
\paragraph{Probability for the independence property}
First, we provide an upper bound for the probability that a pixel has the independence property.

\begin{Le} 
\label{lem:choose-grid}
Let $v\in\Z^2_{s,t}$  and let $\varepsilon<\frac{1}{700}$ be fixed. Then
\[ \Prob{ \neg  {\cal I}_\varepsilon(v)} \le   {95 n^3 \expo{-0.194  n} + \left(19n^2+13 n +4  \right)e^{-n\pi}}. \]
\end{Le}
In the above lemma, we have assumed that $\varepsilon<\frac{1}{700}$
to obtain an upper bound for $\Prob{ \neg  {\cal I}_\varepsilon(v)}$
which is independent of $\varepsilon$. 

\begin{prooft}
Let $N_\varepsilon(v)$ be the number of Delaunay triangles in $\DT({X_n})$ such that
the associated circumdisk intersects simultaneously  
$C^\varepsilon(v)$ and the {complement} of $C_2(v)$. 
If the event ${\cal I}_\varepsilon(v)$ does not hold, then $N_\varepsilon(v) \ge 1$ (here we have used the fact that ${\cal  I}_\varepsilon(v)$ is $\sigma(X_n \cap C(v))$ measurable).
It follows from the Markov's inequality that
\[
\Prob{ \neg  {\cal I}_\varepsilon(v)} 
\le
\Prob{  N_\varepsilon(v) \ge 1} 
\le
\Ex{  N_\varepsilon(v) }. 
\]
Besides, 
\begin{eqnarray*}
\Ex{  N_\varepsilon(v) }  &=& 
{\frac{1}{3!}}
\Ex{\sum_{{p_{1\diff 3}\in X_n^3}}
       \indicator{\Delta(p_{1:3})\in \DT(X_n)}
       \indicator{B(p_{1:3})\not\subset C_2(v);\; B(p_{1:3})\cap C^\varepsilon(v)\neq\varnothing } }. 
\end{eqnarray*}
It follows from the Slivnyak-Mecke and the Blaschke-Petkantschin formulas that
\begin{eqnarray}
\Ex{  N_\varepsilon(v) }  &=&  \notag
\tfrac{n^3}{3!}\int_{(\R^2)^3} \!\!\! \Prob{B(z,r)\cap X_n=\varnothing}
       \indicator{B(z,r)\not\subset C_2(v);\; B(z,r)\cap C^\varepsilon(v)\neq\varnothing }
            dp_{1:3}
\\ &=& \notag
\tfrac{n^3}{6} \int_{\R_+}\int_{\R^2}\int_{\s^3 }
        \expo{-n\Area{B(z,r)}}
       \indicator{B(z,r)\not\subset C_2(v);\; B(z,r)\cap C^\varepsilon(v)\neq\varnothing }
\\ & & \hspace*{2cm} \notag
\,\times\,
       r^3 
          2\Area{\Delta(u_{1:3})}
        \sigma(du_{1:3}) dz\, dr 
\\  &=&\label{eq:majIepsilon}
\tfrac{n^3}{6} 24 \pi^2  \int_{\R_+}\!\int_{\R^2}\!\!\! 
        \expo{-n\pi r^2}
       \indicator{B(z,r)\not\subset C_2(v); z\in C^\varepsilon(v)\oplus B(0,r) } 
       r^3   dz\, dr
\\ & & \notag
\end{eqnarray}
since  $B(z,r)\cap C^\varepsilon(v)\neq\varnothing$ if and only
  if $z\in C^\varepsilon(v)\oplus B(0,r)$. To deal with 
  $B(z,r)\not\subset C_2(v)$ we consider two cases as follows {(see Figure~\ref{fig:choose-grid})}:

\textbf{Case 1:} if $r\leq 1$, we use the fact that 
 \[     
  \left\{
      \begin{split} & B(z,r)\not\subset C_2(v) \\
        &  z\in  C^\varepsilon(v)\oplus B(0,r) \end{split}  \right.
     \Rightarrow       z\in[-r-\tfrac{1}{2}-\epsilon, r+\tfrac{1}{2}+\epsilon]^2   \setminus [r-1,1-r]^2.
\] 
Notice that we have to choose $r>\frac{1 - 2\epsilon}{4}>0.249$ to ensure that the set of the right-hand side is not empty. Besides, this set has an area smaller than $4(2r+\varepsilon-\tfrac{1}{2})(\tfrac{3}{2}+\varepsilon)<12.016 r$. 

\textbf{Case 2:} if $r> 1$, we use the trivial assertion
 \[     
  \left\{
      \begin{split} & B(z,r)\not\subset C_2(v) \\
        &  z\in  C^\varepsilon(v)\oplus B(0,r) \end{split}  \right.
     \Rightarrow    z\in[-r-\tfrac{1}{2}-\epsilon, r+\tfrac{1}{2}+\epsilon]^2.
\] 
The set of the right-hand side has an area which is lower than $(1+2\varepsilon+2r)^2<3.004r^2$.
\begin{figure}[t]
\begin{center}
  \includegraphics[page=10,width=\textwidth]{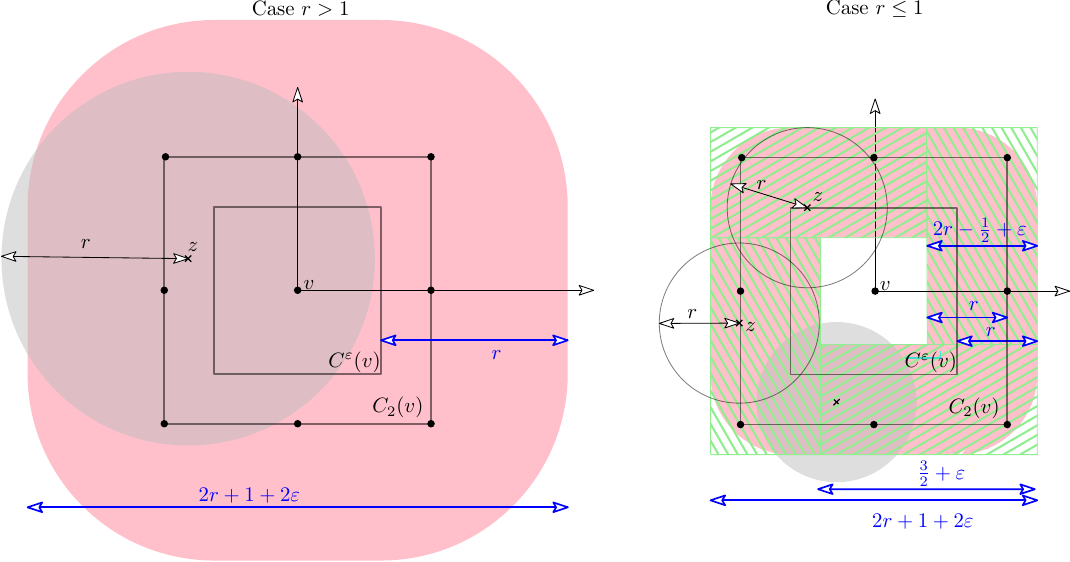}
\end{center}
\caption{
      The squares $C^{\varepsilon}(v)$, $C_2(v)$ (black), the set
        $C^\varepsilon(v)\oplus B(0,r)$ (pink) and the disk $B(z,r)$ (grey).
        \label{fig:choose-grid}
     }
 \end{figure}
 By integrating over $z$, it follows from  \eqref{eq:majIepsilon} that 
\begin{eqnarray*}
\Prob{ \neg  {\cal I}_\varepsilon(v) }
&\leq& 
4\pi^2n^3 
 \pth{
 \int_{0.249}^1 \expo{-n\pi r^2}   12.016 \; r^4 dr 
 + \int_1^\infty  \expo{-n\pi r^2}  3.004 \; r^5 dr   }
 \\ &\leq  &  
 4\pi^2n^3 \pth{\expo{-n\pi 0.249^2}\int_{0.249}^1  12.016 \; r^4 dr  
 + \int_1^\infty  \expo{-n\pi r^2}  3.004 \; r^5 dr   } \\
  &\leq  & 
 95 n^3 \expo{-0.194  n} + \left(19n^2+13 n +4 \right)e^{-n\pi}.
\end{eqnarray*}
\end{prooft}

\paragraph{Probability for the weak horizontality property}
Secondly, we provide an upper bound for the probability that a pixel has a weak horizontality property conditional on the fact that it has the independence property. 

\begin{Le} 
\label{lem:horizontality}
Let $v\in\Z^2_{s,t}$ and $\rho <0.4$. Then there exists $\kappa>0$ such that 
$\tfrac{\kappa}{\kappa-1}\rho<\tfrac{\pi^2}{8}$ and 
\[ \Prob{  {\cal I}_{\varepsilon_\rho}(v)} \Prob{  {\cal H}'_{\varepsilon_\rho,\alpha_{\rho,\kappa},\kappa}(v) \midd  {\cal I}_{\varepsilon_\rho}(v)} 
\le  
31.8\pth{\frac{3}{4}+\frac{\sqrt{\rho}\sqrt{\rho+2}}{2}}^2 \sqrt{\rho n},\]
where   $\varepsilon_\rho:=\sqrt{\rho}\sqrt{2+\rho}$ and 
$\alpha_{\rho,\kappa} := \sqrt{2\tfrac{\kappa}{\kappa-1}\rho}\in [0, \tfrac{\pi}{2}]$ are the same as in Lemma~\ref{lem:alpha}. 
\end{Le}

\begin{prooft}
According to the Markov's inequality, we have:
\begin{align}
 \label{eq:majhorizon1}
\Prob{  {\cal H}'_{\varepsilon_\rho,\alpha_{\rho,\kappa},\kappa}(v) \midd  {\cal I}_{\varepsilon_\rho}(v)} & \notag = \Prob{ L_{\varepsilon_\rho,\alpha_{\rho,\kappa},\kappa}(v)\geq\tfrac{1}{\kappa} \midd  {\cal I}_{\varepsilon_\rho}(v)}\\
& \leq  \notag \kappa \Ex{ L_{\varepsilon_\rho,\alpha_{\rho,\kappa},\kappa}(v)\midd  {\cal    I}_{\varepsilon_\rho}(v)} \\
& =  \frac{\kappa\Ex{ L_{\varepsilon_\rho,\alpha_{\rho,\kappa},\kappa}(v) \indicator{{\cal I}_{\varepsilon_\rho}(v)}}}{\Prob{{\cal  I}_{\varepsilon_\rho}(v)}}.
\end{align}

Now,  recalling that on the event  $ {\cal I}_{\varepsilon_\rho}(v)$,
any triangle, and thus any  edge, intersecting
$C^{\varepsilon_\rho}(v)$ has its vertices inside $C_2(v)$, we have:
\begin{multline*}
\Ex{  L_{\varepsilon_\rho,\alpha_{\rho,\kappa},\kappa}(v)   \indicator{{\cal I}_{\varepsilon_\rho} (v)}}\\
 \leq  \!
\tfrac{1}{2}\!
\Ex{\!\sum_{p_{1\diff 3}\in X_n^3} \!\!\!\!\!
       \indicator{p_{1:3}\in DT(X_n\!)}
       \indicator{\eangle{p_1 p_2}<\alpha_{\rho,\kappa};\, x_{p_1} \leq x_{p_2}}
       \indicator{B(p_{1:3})\subset C_2(v)}
       \indicator{B(p_{1:3})\cap C^{\varepsilon_\rho}(v)\neq \varnothing}\!
              h(p_1p_2)\!
}
\end{multline*}
where the $\tfrac{1}{2}$ factor comes from the fact that each edge is
counted twice in the sum (once for each incident triangle). %% pas si x1<x2
%As in the proofs of Propositions~\ref{prop:origin}, \ref{lem:straight}, and~\ref{lem:upper-bound}, 
We apply the
Slivnyak-Mecke and the Blaschke-Petkantschin formulas. 
This gives:
\begin{multline*}
\Ex{  L_{\varepsilon_\rho,\alpha_{\rho,\kappa},\kappa}(v)   \indicator{{\cal I}_{\varepsilon_\rho} (v)}} \leq \tfrac{1}{2} {n^3} \int_{C_2(v)} \int_{0}^1\int_{\s^3 }
         \expo{-n\Area{B(z,r)}}
   \indicator{\eangle{u_1u_2}<\alpha_{\rho,\kappa} ;\, x_{u_1}\leq  x_{u_2}}
   \\\times     \indicator{B(z,r)\subset C_2(v)}
       \indicator{B(z,r)\cap C^{\varepsilon_\rho}(v)\neq \varnothing}
       \cdot r \,h(u_1u_2) 
       \cdot    r^3 2\Area{\Delta(u_{1:3})}    \sigma(du_{1:3}) dr dz
\end{multline*}
since $B(z,r)\subset C_2(v)$ implies that $z\in C_2(v)$ and $r\leq
1$.
Hence, $\Ex{  L_{\varepsilon_\rho,\alpha_{\rho,\kappa},\kappa}(v)   \indicator{{\cal I}_{\varepsilon_\rho} (v)}} \leq a(n,\rho)\times b(\rho,\kappa)$, where
\begin{equation*}
a(n,\rho) := \tfrac{1}{2} {n^3} \int_{C_2(v)} \int_{0}^1
         \expo{-n\pi r^2}r^4   \indicator{B(z,r)\subset C_2(v)}
       \indicator{B(z,r)\cap C^{\varepsilon_\rho}(v)\neq \varnothing}
           dr dz.
\end{equation*}
\begin{equation*}
b(\rho,\kappa):=\int_{\s^3} 
         \indicator{\eangle{u_1u_2}<\alpha_{\rho,\kappa} ;\, x_{u_1}\leq  x_{u_2}} 2\Area{\Delta(u_{1:3})} h(u_1u_2) \sigma(du_{1:3}).
\end{equation*}

First, we provide an upper bound for $a(n,\rho)$. 
Since $B(z,r)\cap C^{\varepsilon_\rho}(v)\neq \varnothing$ 
we have $z\in
v\oplus[-\frac{3}{4}-\frac{\varepsilon_\rho}{2},\frac{3}{4}+\frac{\varepsilon_\rho}{2}]^2$. 
Dividing the square 
$v\oplus[-\frac{3}{4}-\frac{\varepsilon_\rho}{2},\frac{3}{4}+\frac{\varepsilon_\rho}{2}]^2$  
into four  quadrants of equal size,  it follows that
\begin{figure}[t]
\begin{center}
  \includegraphics[page=11,width=\textwidth]{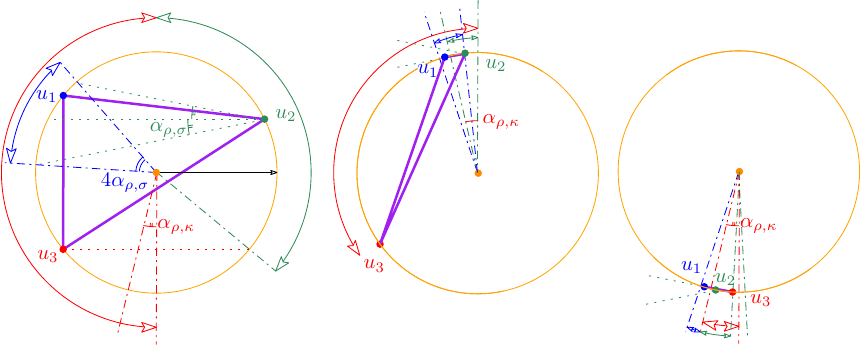}
\end{center}
\caption{
      Domains of integration for $u(\beta_{1:3})$.
        \label{fig:domain-alpha}}
 \end{figure}

\begin{align}
\label{eq:majhorizon2}
a(n,\rho) & \notag \leq 
 \tfrac{1}{2} 4   n^3 \int_{0}^{\frac{3}{4}+\frac{\varepsilon_\rho}{2}}
\int_0^{\frac{3}{4}+\frac{\varepsilon_\rho}{2}}
%\int_{\max(0,x_z-\tfrac{1}{2}-\varepsilon_\rho,y_z-\tfrac{1}{2}-\varepsilon_\rho)}^{\min(1-x_z,1-y_z)}\hspace*{-1cm}
\int_0^\infty
\expo{-n\pi r^2}r^4 dr dy_z dx_z\\
& = 2n^3\pth{\frac{3}{4}+\frac{\varepsilon_\rho}{2}}^2 \,\frac{3}{8\pi^2n^{\frac{5}{2}}}. 
\end{align}

Secondly, to provide an upper bound for $b(\rho,\kappa)$,  we write  
\[u_{1:3}:=u(\beta_{1:3}):=(u(\beta_1),u(\beta_2),u(\beta_3))\] with
$u(\beta_i)=(\cos\beta_i, \sin\beta_i)$. Up to the line symmetry
w.r.t. the $y$-axis, we impose that
$\beta_3\in[\tfrac{\pi}{2},\tfrac{3\pi}{2}]$. Up to the line symmetry
w.r.t. the $x$-axis, we also impose that 
$y_{u(\beta_2)}\geq y_{u_{\beta_3}}$, i.e. 
$\beta_2\in [\pi-\beta_3,\beta_3]$. 
Besides,  $u_1\in C(0,1)\cap C(\beta_2)$, where $C(0,1)$ is the unit circle and where $C(\beta_2)$ is the half-cone generated on the
left of $u(\beta_2)$ by the $x$-axis, with
vertex $u(\beta_2)$ and with angle $2\alpha_{\rho,\kappa}$, i.e.
\[C(\beta_2):=u(\beta_2)+\{(r\cos\gamma, r\sin\gamma): r\geq 0, \gamma\in [\pi-\alpha_{\rho,\kappa}, \pi+\alpha_{\rho,\kappa}]\}. \]
We
discuss four cases by splitting  the domain of integration 
$ [\pi-\beta_3,\beta_3]$ 
of $\beta_2$ as follows  (see Figure~\ref{fig:domain-alpha}):
\begin{enumerate}
\item 
If $\beta_2\in[\pi-\beta_3,\tfrac{\pi}{2}]$, we have $\beta_1\in[\pi-\beta_2-2\alpha_{\rho,\kappa},\pi-\beta_2+2\alpha_{\rho,\kappa}]$ ( Figure~\ref{fig:domain-alpha}-left).
\item 
If $\beta_2\in[\tfrac{\pi}{2},\tfrac{\pi}{2}+\alpha_{\rho,\kappa}]$, 
 the length of $C(0,1)\cap C(\beta_2)$  is maximal when
$\beta_2=\frac{\pi}{2}$, so that $\beta_1\in
[\tfrac{\pi}{2},\tfrac{\pi}{2}+2\alpha_{\rho,\kappa}]$. Besides, 
the area of the triangle $\Delta(u(\beta_{1:3}))$ is less than $2\alpha_{\rho,\kappa}$
and $h(u(\beta_1)u(\beta_2))$ is also  less than $2\alpha_{\rho,\kappa}$ 
(Figure~\ref{fig:domain-alpha}-center). 
\item 
If $\beta_2\in[\tfrac{\pi}{2}+\alpha_{\rho,\kappa},\tfrac{3\pi}{2}-\alpha_{\rho,\kappa}]$, we have $C(0,1)\cap C(\beta_2)=\varnothing$.
\item 
If $\beta_2\in[\tfrac{3\pi}{2}-\alpha_{\rho,\kappa},\beta_3]$, with $\beta_3>\tfrac{3\pi}{2}-\alpha_{\rho,\kappa}$, the length of $C(0,1)\cap C(\beta_2)$  is maximal when
$\beta_2=\frac{3\pi}{2}$, so that $\beta_1\in
[\tfrac{3\pi}{2}-2\alpha_{\rho,\kappa},\tfrac{\pi}{2}]$.
Besides, the area of the triangle  $\Delta(u(\beta_{1:3}))$ is less than $2\alpha_{\rho,\kappa}$
and $h(u(\beta_1)u(\beta_2))$ is also less than $2\alpha_{\rho,\kappa}$ (Figure~\ref{fig:domain-alpha}-left).
\end{enumerate}
It follows that
\begin{multline*}
b(\rho,\kappa) 
\leq  
4
\int_{\frac{\pi}{2}}^{\frac{3\pi}{2}}\int_{\pi-\beta_3}^{\frac{\pi}{2}}
\int_{\pi-\beta_2-2\alpha_{\rho,\kappa}}^{\pi-\beta_2+2\alpha_{\rho,\kappa}}
2 \Area{\Delta(u(\beta_{1:3}))} h({u(\beta_1)u(\beta_2)})
         d\beta_1d\beta_2d\beta_3 + 16\pi\alpha_{\rho,\kappa}^4. 
\end{multline*}
Moreover, we know that
\[2 \Area{\Delta(u(\beta_{1:3}))} = \betamatrix \quad \text{and} \quad h({u(\beta_1)u(\beta_2)}) = \cos \beta_1 \!-\! \cos \beta_2.\]
This gives 
\begin{equation*}
\begin{split}
b(\rho,\kappa) & \leq   \tfrac{256}{9}  \cos^3  \alpha_{\rho,\kappa} \sin \alpha_{\rho,\kappa}  +\tfrac {128}{3}\cos\alpha_{\rho,\kappa} \sin\alpha_{\rho,\kappa}
+\tfrac{128}{3}\alpha_{\rho,\kappa}+ 16\pi\alpha_{\rho,\kappa}^4\\ 
& \leq  \tfrac{1024}{9}\alpha_{\rho,\kappa}.
 \end{split}
\end{equation*}

This together with \eqref{eq:majhorizon1}, \eqref{eq:majhorizon2}  and the fact that 
$\Ex{  L_{\varepsilon_\rho,\alpha_{\rho,\kappa},\kappa}(v)
  \indicator{{\cal I}_{\varepsilon_\rho} (v)}} \leq a(n,\rho)\times b(\rho,\kappa)$ 
with 
$\varepsilon_\rho=\sqrt{\rho}\sqrt{\rho+2}$ and $\alpha_{\rho,\kappa}=\sqrt{2\frac{\kappa}{\kappa-1}\rho}$, implies that 
\begin{equation*}
 \begin{split}
 \Prob{  {\cal I}_{\varepsilon_\rho}(v)} \Prob{ {\cal H}'_{\varepsilon_\rho,\alpha_{\rho,\kappa},\kappa}(v) \midd  {\cal I}_{\varepsilon_\rho}(v) } 
  & \le \kappa \cdot  2 n^3  \pth{\frac{3}{4}+\frac{\sqrt{\rho}\sqrt{\rho+2}}{2}}^2 \,\frac{3}{8\pi^2n^{\frac{5}{2}}}\cdot \tfrac{1024}{9}\cdot  \sqrt{2\frac{\kappa}{\kappa-1}\rho}\\
  & = 2 \cdot  \frac{3}{8\pi^2}\tfrac{1024}{9} \sqrt{2}  \pth{\frac{3}{4}+\frac{\sqrt{\rho}\sqrt{\rho+2}}{2}}^2 \sqrt{\tfrac{\kappa^3}{\kappa-1}}\sqrt{\rho n}\\
   & \le   31.76 \pth{\frac{3}{4}+\frac{\sqrt{\rho}\sqrt{\rho+2}}{2}}^2 \sqrt{\rho n}.
   \end{split}
\end{equation*}
In the last line, we have taken $\kappa=\tfrac{3}{2}$ since it is the value of $\kappa$ which minimizes
$\tfrac{\kappa^3}{\kappa-1}$. The condition $\tfrac{\kappa}{\kappa-1}\rho<\tfrac{\pi^2}{8}$
is satisfied since $\rho<\frac{\pi^2}{24}\simeq 0.4$ by assumption.
\end{prooft}

\begin{prooft}{ Proposition~\ref{lem:good-proba}}
This is a direct consequence of  Lemmas~\ref {lem:alpha}~,\ref{lem:choose-grid} and~\ref{lem:horizontality} and the fact that 
\begin{align*}
\Prob{  \neg{\cal I}_{\varepsilon_\rho}(v) \vee {\cal H}_\rho (v) }
&\le 
\Prob{ {\neg {\cal I}_{\varepsilon_\rho} (v)} }
+ \Prob{ {\cal I}_{\varepsilon_\rho} (v)}
\Prob{{\cal H}'_{\varepsilon_\rho,\alpha_{\rho,\kappa},\kappa}(v)\midd {\cal I}_{\varepsilon_\rho} (v)}
\end{align*}
The assumption $\rho<4\cdot 10^{-6}$  ensures that
$\varepsilon_\rho<\frac{1}{700}$ in Lemma~\ref{lem:choose-grid}
and $\rho<0.4$ in Lemma~\ref{lem:horizontality}.
\end{prooft}

\subsection{Wrap-up\label{sec:wrapup}}

\begin{prooft}{Theorem~\ref{th:lower}}
 It follows from Lemma~\ref{lem:length-animal} that 
\begin{align*}  
 \Prob{ \ell(SP_X)  \le  k + \rho  \pth{k - 16  k \sqrt{p} }}
& \leq 
 \Prob{ \max_{c\in Colors} \sharp_{{\cal H}} (\animal_{(c)}(SP_X))  \ge  4k \sqrt{p} }
\\ &\leq
\sum_{c\in Colors}  \Prob{ \sharp_{{\cal H}} (\animal_{(c)}(SP_X))  \ge 4 k \sqrt{p} }.
\end{align*}
To bound the right-hand side, we apply
Proposition~\ref{lem:percolation} to the family of events
$\mathcal{Y}:=(Y_v)_{v\in\Z^2}$, where $Y_v:={\cal H}_\rho(v) \vee \neg {\cal I}_{\varepsilon_\rho}(v)$. 
Notice that for any color $c$, the events $(Y_v)_{v\in \Z^{2,
    (c)}_{s,t}}$ are independent: this comes from  the fact that the
event 
$ {\cal H}_\rho(v) \vee \neg{\cal I}_{\varepsilon_\rho} (v) $ is $\sigma(X_n\cap C_2(v))$ 
measurable and the fact that $C_2(v)\cap C_2(w)=\varnothing$ when
$v\neq w\in\Z^{2,(c)}_{s,t}$.   
Moreover, for any $v\in\Z^{2,(c)}_{s,t}$, the probability
 $p:=\Prob{ {\cal H}_\rho(v) \vee  \neg {\cal  I}_{\varepsilon_\rho}(v)}$ does not depend on $v$. 
Besides, to apply Proposition~\ref{lem:percolation}, we have to choose
$\rho$ and $n$ in such a way that $p\leq 0.01$. Provided that such a
condition is satisfied and that $\rho\leq 4\cdot 10^{-6}$, it follows
from Proposition~\ref{lem:percolation} and
Proposition~\ref{lem:good-proba} that  
\[ \Prob{ \ell(SP_X)  \le  k + \rho  \pth{1 - 16   \sqrt{P(\rho, n)} }k } = O(k^{-1/2}).\]

\noindent\begin{minipage}{0.5\textwidth}
To minimize $\rho  \pth{1 - 16   \sqrt{P(\rho, n)} }$, we choose
$\rho=1.25\cdot10^{-10}$ and $n=153$. 
Then we have $\rho <4\cdot 10^{-6}$ and 
$p\leq P(\rho,n)\simeq 0.0025<0.01$, 
so that the assumptions of 
Propositions~\ref{lem:percolation} and \ref{lem:good-proba} 
are satisfied. This concludes the proof of Theorem~\ref{th:lower} since 
\[\pth{1 - 16   \sqrt{P(\rho, n)} } \geq 2.47\cdot 10^{-11}.\]
\end{minipage}\quad\begin{minipage}{0.4\textwidth}
\includegraphics[page=12,width=\textwidth]{Figures}
\end{minipage}
\end{prooft}

\section{Simulations\label{s:simulation}}
We give below experimental values,  which depend on the intensity $n$, for
the expected length and the number of edges of the paths considered in
this paper. 
We also experiment a third path, described below:\\
\textit{Greedy constructed Path $GP_X$:} we define such a path, starting 
in $s$, by induction. 
Let $(T_i)_{0\leq i\leq k}$  be the family 
of triangles in $\DT(X)$  whose interior  intersects $[s,t]$ and ordered from the 
left to the right. 
Let $w$ be the last vertex inserted in the path and let $i$ be the largest value such that  $w\in T_i$. 
Then the edge of $T_i$ incident to $w$ which minimizes the angle with respect to 
the $x$-axis is added to the path (see Figure~\ref{fig:paths}).

\begin{figure}[t]
\noindent
     \begin{center}
         \includegraphics[page=14,width=0.5\textwidth]{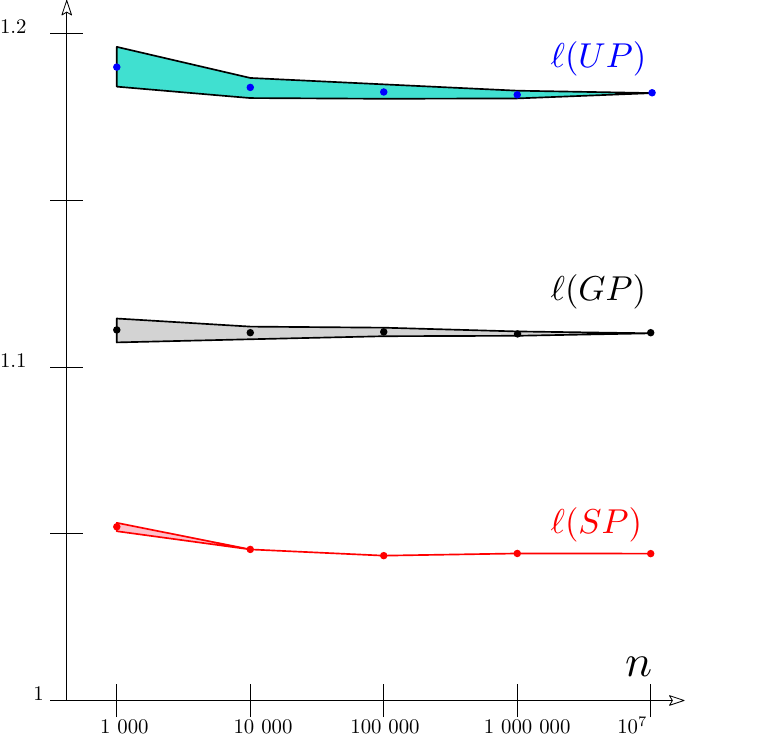}
     \end{center}
   \caption{
       Experimental evaluations of $\Ex{\length{SP_X}}$,
       $\Ex{\length{GP_X}}$, and $\Ex{\length{UP_X}}$
      for various intensities with  $\|s-t\|=1$. The width of the drawing is the
      standard deviation.
        \label{fig:experiments}
     }
\end{figure}

The code for these simulations is written with {CGAL}~\cite{cgal}
and is available in~\cite{hal-version}.
Figure~\ref{fig:experiments} and Table \ref{Table:experiments}  give estimates of the expectations and of the standard deviations for the lengths and for the cardinalities of the paths $UP_X$, $GP_X$ and $SP_X$. These estimates are based on 100 trials of Poisson point processes with intensity $n=$10 millions.
In particular, Propositions~\ref{lem:straight} and~\ref{lem:upper-bound} are confirmed by
our experiments. The correct value for $\Ex{\length{SP_X}}$ is about 1.04
and the path $GP_X$, using simple trick to improve on $UP_X$,
gives a path significantly shorter than  $UP_X$. We also give experimental values for the straigth walk $SW_X$.

     \begin{table}[h!]
\begin{center}
       \begin{tabular}{|c|l|l|c|c|}\hline 
\rule[-4mm]{0mm}{11mm}         path $P$ & $\Ex{\length{P}}$ & $\sigma(\length{P})$ & $\displaystyle\frac{\Ex{\size{P}}}{\sqrt{n}}$ & $\displaystyle\frac{\sigma(\size{P})}{\sqrt{n}}$\\\hline \hline 
 \rule{0mm}{4mm}        $SW_X$ &            &              & 2.1602 & 0.0154 \\\hline 
  \rule{0mm}{4mm}        $UP_X$ & 1.1826 & 0.0053 & 1.0804 & 0.0102 \\\hline 
   \rule{0mm}{4mm}       $GP_X$ & 1.1074 & 0.0036 & 1.0062 & 0.0084 \\\hline 
  \rule{0mm}{4mm}        $SP_X$ & 1.0401 & 0.0004 &  0.9249 & 0.0080 \\\hline 
       \end{tabular}\\
            \end{center}
            \caption{Experimental values for $SW_X$, $UP_X$, $GP_X$ and $SP_X$        \label{Table:experiments} }
       \end{table}
       
\begin{figure}[p]
     \begin{center}
  \includegraphics[page=13, scale=0.6]{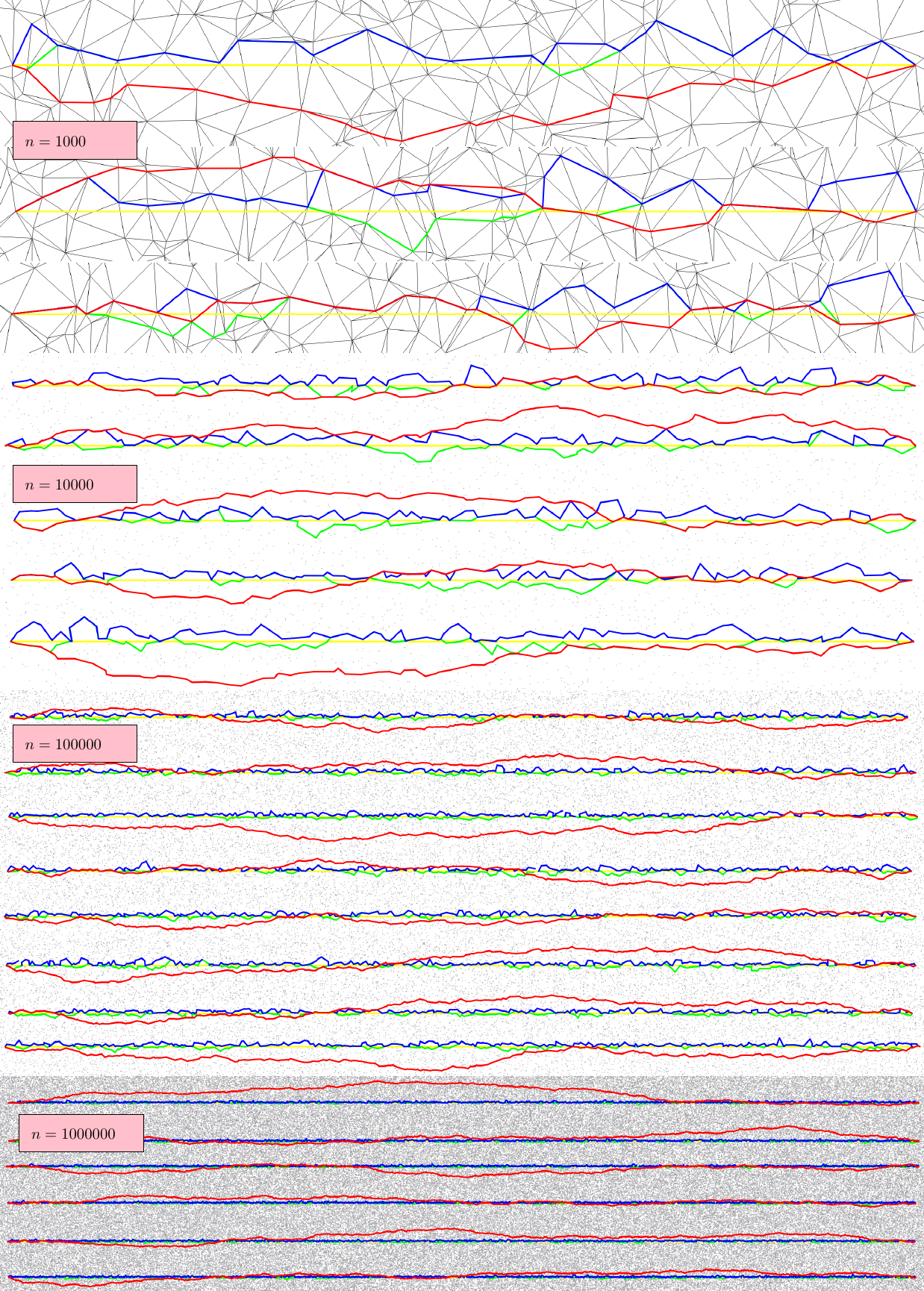}
     \end{center}
     \caption{
       Several  trials with intensity from $n=1000$ to $n=1$   million.
       The paths $UP_X$, $GP_X$,
       and $SP_X$ are colored in blue, green, and red respectively and
       the line segment  $[s,t]$ is colored in yellow. 
        \label{fig:trials}
     }
 \end{figure}

Several illustrations of the paths $UP_X$, $GP_X$ and $SP_X$ are depicted in Figure~\ref{fig:trials} for various intensities.

\section{Conclusion}
\label{s:conclusion}
We have provided a non-trivial lower bound for the {inferior} limit of the expectation of
  the stretch factor between 
two points $s$ and
  $t$ in $Del(X_n\cup\{s,t\})$ where $X_n$ is an homogeneous Poisson point process
when the intensity goes to infinity.
The main difficulty is that  we have absolutely no information on
the exact location of the shortest path.
Our lower bound for $\rho$ is far from being tight since the experimental
value is much larger.
Although several constants in our  proof can be improved a little bit,
this scheme of proof can only give lower bounds which are far from optimal. Indeed, 
the first point where our evaluation is quite crude is the approximation
of the Horizontality Property by the Weak Horizontality Property. Actually, a better upper bound for 
$\Prob{{\cal H}_\rho(v)}$ should widely improve $\rho$ 
and an explicit lower bound for 
$\Prob{{\cal H}_\rho(v)}$ yields to a more precise evaluation of the
convergence.
Another point where we widely under-evaluate $\rho$ comes from the fact that we use approximation by animals. Actually, a bad situation (when the shortest path is very short) corresponds to an animal with many 
pixels with a strong horizontality property. However, the converse is clearly not true. Thus we believe that the
proof of a tight constant necessitates other techniques. 

Another question raised by the paper is to prove that the limit of the
expectation of the stretch factor exists. One of the underlying
difficulties is to show that the function $k\mapsto
\Ex{\length{SP_{X}}}$, with $X=X_n\cup\{(0,0), (k,0)\}$ and
$k\in\R_+$, is continuous and subadditive. It appears that these
properties are not trivial. 
Indeed, conditioning by the fact that an intermediate point in $[s,t]$
belongs to the shortest path can increase or decrease the length of
the shortest path in distribution. We hope, in a future work, to deal
with these properties and to prove the existence of the limit.

\bigskip

\textbf{Acknowledgement}
The authors thanks David Coeurjolly for pointing us ref~\cite{gerard2015tight}
for Lemma~\ref{lem:size-animal}, and Louis Noizet for discussions
about the definition of $UP_X$.

%%%%%%%%%%%%%% INSERT bbl FILE FOR FINAL SUBMISSION

%\bibliographystyle{amsplain}
%\bibliography{biblio}

\def\doi#1{\href{http://dx.doi.org/#1}{doi:#1}}
\providecommand{\bysame}{\leavevmode\hbox to3em{\hrulefill}\thinspace}
\providecommand{\MR}{\relax\ifhmode\unskip\space\fi MR }
\providecommand{\MRhref}[2]{%
  \href{http://www.ams.org/mathscinet-getitem?mr=#1}{#2}
}
\providecommand{\href}[2]{#2}

\newpage
\appendix

\section{Proof for the variance of $UP_{X}$\label{ap:variance}}

\begin{prooft}[Proof of Proposition~\ref{prop:up-variance}]
{According to Proposition~\ref{lem:upper-bound} and the fact that $ \Ex{\ell(UP_{X})} = \Ex{L_{X_n}} +  O\pth{n^{-\tfrac{1}{2}}}$, where $X=X_n\cup\{s,t\}$ and $s=(0,0)$, $t=(1,0)$, it is enough to show that $\Ex{L_{X_n}^2}=\left(\frac{35}{3\pi^2}\right)^2 + O\pth{n^{-\tfrac{1}{2}}}$.}

 Thanks
to \eqref{def:Lchi}, we have $\Ex{L_{X_n}^2}  = \Ex{\left(\sum_{p_{1\diff 3}\in
    X_n^3} \ell_{X_n}(p_{1:3}) \right)^2}$. To estimate the right-hand
side, we will develop the terms inside the sum by discussing the
number of common vertices  between two
triangles. To do it, we first give some notation. For each $0\leq
  l\leq 2$, we denote by $\mathcal{S}_l\subset \{1:3+l\}^3$ the set: 
\[\mathcal{S}_l:=\{(i_1,i_2,i_3)\in \{1:3+l\}^3_{\neq} : \size{\{i_1,i_2,i_3\}\cap \{1,2,3\}} = 3-l   \}.\]
Given a  $(3+l)$-tuple of points $p_{1:3+l}$, we can associate a
family of couple of triangles as follows. The first triangle is
$\Delta(p_{1:3})$ and the second one is any triangle of the form
$\Delta(p_{i_1}p_{i_2},p_{i_3})$ such that $(i_1,i_2,i_3)\in
\mathcal{S}_l$. In particular, when $l=0$, we have
$\mathcal{S}_0=\{1,2,3\}^3_{\neq}$ and $\Delta(p_{i_1},p_{i_2},
p_{i_3})=\Delta(p_{1:3})$.  Notice that the term $l$ is the number of
vertices of the triangle $\Delta(p_{i_1},p_{i_2}, p_{i_3})$ which are
not vertices of the triangle $\Delta(p_{1:3})$. By developing the sum
associated with $L_{X_n}$, we get
\begin{equation*}
L_{X_n}^2  = \sum_{p_{1\diff 6} \in X_n^6}l_{X_n}(p_{1:3})l_{X_n}(p_{4:6})
+ \sum_{l=0}^2  \sum_{(i_1,i_2,i_3)\in\mathcal{S}_l} 
             \sum_{p_{1\diff 3+l} \in X_n^{3+l}}
                         l_{X_n}(p_{1:3})l_{X_n}(p_{i_1}, p_{i_2}, p_{i_3}).
\end{equation*}

As a first step, we show that 
\[
\Ex{\sum_{p_{1\diff 6} \in X_n^6
  }l_{X_n}(p_{1:3})l_{X_n}(p_{4:6})} =
\left(\frac{35}{3\pi^2}\right)^2+O\pth{n^{-\tfrac{1}{2}}}.
\]
 In the same spirit as in the proof of Proposition~\ref{lem:upper-bound}, we apply the Slivnyak-Mecke and Blaschke-Petkantschin fomulas. It follows that
\begin{multline*}  
\Ex{\sum_{(p_{1\diff 6}) \in X_n^6   }l_{X_n}(p_{1:3})l_{X_n}(p_{4:6})} = n^6\cdot \int_{(\R^2)^2} \int_{\R_+^2}\int_{(\s^3)^2} e^{-n \Area{B(z,r)\cup B(z',r')}}\\
\begin{split}
& \times \indicator{z+ru_{1:3}\in E^+}\indicator{z'+ru'_{1:3}\in E^+} r^4r'^4\|u_2-u_1\|\cdot \|u'_2-u'_1\|
\\ 
& \times \Area{\Delta(u_{1:3})}\Area{\Delta(u'_{1:3})}\sigma(du_{1:3})\sigma(du'_{1:3})drdr'dzdz'.\end{split} 
\end{multline*}
Rewriting the exponentiel term in the integrand as  \[e^{-n \Area{B(z,r)\cup B(z',r')}} = e^{-n\pi (r^2+r'^2)} + \left(e^{-n \Area{B(z,r)\cup B(z',r')}}-e^{-n\pi (r^2+r'^2)}\right)\] and applying Fubini's theorem, we deduce from \eqref{eq:Upintegral} that 
 \begin{equation}\label{eq:vardisjoint}
\Ex{\sum_{(p_{1\diff 6}) \in X_n^6   }l_{X_n}(p_{1:3})l_{X_n}(p_{4:6})} = (\Ex{L_{X_n}})^2 + r_n,
\end{equation} 
 where 
 \begin{multline*}
 r_n:=  2n^6\cdot \int_{(\R^2)^2} \int_{\R_+^2}\int_{(\s^3)^2}      \left(e^{-n \Area{B(z,r)\cup B(z',r')}}-e^{-n\pi (r^2+r'^2)}\right)\\
 \begin{split} & \times\indicator{z+ru_{1:3}\in E^+}\indicator{z'+ru'_{1:3}\in E^+}
 r^4r'^4 \|u_2-u_1\|\cdot \|u'_2-u'_1\|\Area{\Delta(u_{1:3})}\\
&\times \Area{\Delta(u'_{1:3})}\indicator{B(z,r)\cap
  B(z',r')\neq\varnothing}\indicator{r'\leq r}
\sigma(du_{1:3})\sigma(du'_{1:3})dzdz'drdr'.\end{split}  \end{multline*}
Notice that in the integrand of the rest $r_n$, we have added two
indicator functions. The first one deals with the event 
$\{B(z,r)\cap B(z',r')\neq\varnothing\}$: 
this comes from the fact that, on the 
{complement} of this event, the integrand appearing in $r_n$ equals
0 since $\Area{B(z,r)\cup B(z',r')}=\pi (r^2+r'^2)$. The second one concerns the event $r'\leq r$: by symmetry, this explains why we
have added a constant 2 on the left of term $n^6$. Besides, for the
first term of the right-hand side in \eqref{eq:vardisjoint}, we have
$(\Ex{L_{X_n}})^2 =
\left(\frac{35}{3\pi^2}\right)^2+O\pth{n^{-\tfrac{1}{2}}}$ according to
Proposition~\ref{lem:upper-bound}.  
To deal with second term of the right-hand side in
\eqref{eq:vardisjoint}, we bound the terms appearing in the integrand
as follows:
first, we use the fact that 
$e^{-n \Area{B(z,r)\cup    B(z',r')}}-e^{-n\pi (r^2+r'^2)}\leq e^{-n\pi r^2}$ with $r'\leq r$. 
Secondly, for each $v_{1:3}\in \s^3$, we bound
$\|v_2-v_1\|\Area{\Delta(v_{1:3})}$ by a constant, by taking
successively $v_{1:3}=u_{1:3}$ and $v_{1:3}=u'_{1:3}$.  
Thirdly, as in
the proof of Proposition~\ref{lem:upper-bound}, we notice that 
$z\in [0,1]\times [-r,r]$ since $z+ru_{1:3}\in E^+$. 
Finally, we also use
the fact that $\indicator{B(z,r)\cap B(z',r')\neq \varnothing}\leq
\indicator{z'\in B(z,2r)}$. It follows that 
\[r_n\leq cn^6\cdot \int_{(\R^2)^2} \int_{\R_+^2}e^{-n\pi
  r^2}r^4r'^4\indicator{z\in [0,1]\times [-r,r]}\indicator{z'\in
  B(z,2r)}\indicator{r'\leq r}dzdz'drdr'. \]
Integrating successively over $r'\leq r$, $z'\in B(z,2r)$ and $z\in [0,1]\times [-r,r]$, we get 
\[r_n\leq c\cdot n^6\int_{\R_+}r^{12}e^{-n\pi r^2}dr = O\pth{n^{-\tfrac{1}{2}}}.\]

As a second step, we show that 
\[
\Ex{\sum_{p_{1\diff 3+l} \in    X_n^{3+l}}l_{X_n}(p_{1:3})l_{X_n}(p_{i_1}, p_{i_2}, p_{i_3})} =
O\pth{n^{-\tfrac{1}{2}}}
\]
for any $0\leq l\leq 2$ and $(i_1,i_2,i_3)\in\mathcal{S}_l$.
In this case, the triangles $\Delta(p_{1:3})$ and $\Delta(p_{i_1}, p_{i_2}, p_{i_3})$ share at less one vertex in common. By symmetry, we have
\begin{multline*}
\sum_{p_{1\diff 3+l} \in X_n^{3+l}}l_{X_n}(p_{1:3})l_{X_n}(p_{i_1}, p_{i_2}, p_{i_3})\\  =
2\cdot  \sum_{p_{1\diff 3+l} \in X_n^{3+l}}l_{X_n}(p_{1:3})
            l_{X_n}(p_{i_1}, p_{i_2}, p_{i_3})\indicator{R(p_{i_1}, p_{i_2}, p_{i_3})\leq R(p_{1:3}))}
\end{multline*}
Assuming that $R(p_{i_1}, p_{i_2}, p_{i_3})\leq R(p_{1:3})$, we obtain
that  $  l_{X_n}(p_{i_1}, p_{i_2}, p_{i_3})\leq 2\cdot
R(p_{1:3})$. Moreover, we have $p_{i_1},p_{i_2},p_{i_3}\in 3
B(p_{1:3})$ where $3 B(p_{1:3})$ is the ball concentric with $B(p_{1:3})$ and radius three
  times bigger. It follows that
\begin{equation*}
\sum_{p_{1\diff 3+l} \in X_n^{3+l}}l_{X_n}(p_{1:3})l_{X_n}(p_{i_1}, p_{i_2}, p_{i_3}) 
 \leq c\cdot \sum_{p_{1\diff 3+l} \in X_n^{3+l}}l_{X_n}(p_{1:3})
            R(p_{1:3})\indicator{p_{i_1},p_{i_2},p_{i_3}\in 3 B(p_{1:3})}.
\end{equation*}
In the above equation, we have bounded the indicator function $\indicator{R(p_{i_1}, p_{i_2}, p_{i_3})\leq R(p_{1:3}))}$ by 1. In the same spirit as above, we first apply the Slivnyak-Mecke formula. This gives 
\begin{multline*}
\Ex{\sum_{p_{1\diff :3+l} \in    X_n^{3+l}}
             l_{X_n}(p_{1:3})l_{X_n}(p_{i_1}, p_{i_2}, p_{i_3})}\\
  \begin{split}    &       \leq   c\cdot n^{3+l}\int_{{(\R^2)}^3}\left(\int_{(\R^2)^l}\indicator{q_{1:l}\in (3B(p_{1:3}))^l}dq_{1:l}    \right) \Ex{l_{X_n\cup\{p_{1:3}\}}(p_{1:3})}R(p_{1:3})dp_{1:3}\\
  & = c'\cdot n^{3+l}\int_{{(\R^2)}^3}\Ex{l_{X_n\cup\{p_{1:3}\}}(p_{1:3})}(R(p_{1:3}))^{1+2l}dp_{1:3}.
\end{split}   
 \end{multline*}
 It follows from the Blaschke-Petkantschin formula that
 \begin{multline*}
\Ex{\sum_{p_{1\diff :3+l} \in    X_n^{3+l}}
             l_{X_n}(p_{1:3})l_{X_n}(p_{i_1}, p_{i_2}, p_{i_3})}\leq 
 c'\cdot n^{3+l}  \int_{\R_+}\int_{\R^2}\int_{\s^3} e^{-n\pi r^2}
  \|u_2-u_1\| r\\\cdot  \indicator{z+ru_{1:3}\in E^+}  \cdot r^{1+2l}
  \cdot   2 \Area{\Delta(u_{1:3})} r^3 \sigma(du_{1:3})dzdr.              
\end{multline*}            
Since $z+ru_{1:3}\in E^+$ implies that $z\in [0,1]\times [-r,r]$, we
obtain by integrating over $z$ and $u_{1:3}$ 
and using Equations~\eqref{int:pi-r-6},~\eqref{int:pi-r-8}, and~\eqref{int:pi-r-10}
that
\[\Ex{\sum_{p_{1\diff :3+l} \in    X_n^{3+l}}
             l_{X_n}(p_{1:3})l_{X_n}(p_{i_1}, p_{i_2}, p_{i_3})} \leq c''\cdot n^{3+l}\int_{\R_+}e^{-n\pi r^2}r^{6+2l}dr = O\left(n^{-\tfrac{1}{2}} \right).\]
 
\end{prooft}

\section{Integrals\label{ap:integrals}}
In this section, we provide values for several integrals which are often used in the
paper.
These integrals can be computed by using tedious classical
computations. These computations are done with Maple (a Maple sheet is
available in~\cite{hal-version}). 
\begin{eqnarray}
% \label{int:pi-r-3}
%\int_0^\infty \expo{-n\pi r^2}r^3 dr  & = &   \frac{1}{2\pi^2 n^2 }
% \\ \label{int:pi-r-odd}
%\int_0^\infty \expo{-n\pi r^2}r^{2i+1} dr  & = &   \frac{i!}{2\pi^{i+1} n^{i+1} }
% \\ \label{int:pi-r-even}
%\int_0^\infty \expo{-n\pi r^2}r^{2i} dr  & = &   O\pth{\frac{1}{n^{i} \sqrt{n} }}
%\\ 
\label{int:pi-r-4}
\int_0^\infty \expo{-n\pi r^2}r^4 dr  & = &   \frac{3}{8\pi^2 n^{2} \sqrt{n}}
\\  \label{int:pi-r-5}
\int_0^\infty \expo{-n\pi r^2}r^5 dr  & = &   \frac{1}{\pi^3 n^3} 
\\  \label{int:pi-r-5-1-infinity}
\int_1^\infty \expo{-n\pi r^2}r^5 dr  & = &   e^{-n\pi}\pth{\frac{1}{2\pi n} +\frac{1}{\pi^2 n^2}+ \frac{1}{\pi^3 n^3} }
\\  \label{int:pi-r-6}
 \int_0^\infty \expo{-n\pi r^2}r^6 dr  & = &   \frac{15}{16\pi^3 n^3\sqrt{n}} 
%\\  %\label{int:pi-r-7}
% \int_0^\infty \expo{-n\pi r^2}r^7 dr  & = &   \frac{3}{\pi^4 n^4 }
\\  \label{int:pi-r-8}
 \int_0^\infty \expo{-n\pi r^2}r^8 dr  & = &   \frac{105}{32\pi^4 n^4\sqrt{n}} 
%\\  \label{int:pi-r-9}
% \int_0^\infty \expo{-n\pi r^2}r^9 dr  & = &   \frac{12}{\pi^5 n^5 }
\\  \label{int:pi-r-10}
 \int_0^\infty \expo{-n\pi r^2}r^{10} dr  & = &   \frac{945}{64\pi^5 n^5\sqrt{n}} 
%\\  \label{int:pi-r11}
% \int_0^\infty \expo{-n\pi r^2}r^{11} dr  & = &   \frac{60}{\pi^6 n^6 }
\\  \label{int:pi-r-12}
 \int_0^\infty \expo{-n\pi r^2}r^{12} dr  & = &   \frac{10395}{128\pi^6 n^6\sqrt{n}} 
\end{eqnarray}

\begin{align}
\notag & \int_{[0,2\pi)^3} \left| \det\betamatrix   \right| d\beta_{1:3}
\\ &\qquad = 
 3!  \int_0^{2\pi} \int_0^{\beta_3} \int_0^{\beta_2}  \det\betamatrix    d\beta_1d\beta_2d\beta_3
  =  24\pi^2\qquad\qquad\quad
\label{int:aire}  
\\ \notag  & \int_{-1}^{1} \int_{\pi+\smallasin h}^{2\pi-\smallasin h}
     \int_{-\smallasin h}^{\pi+\smallasin h}
      \int_{-\smallasin h}^{\pi+\smallasin }
\!\!  \left | \det \betamatrix \right|
      d\beta_1  d\beta_2  d\beta_3 dh
\\ \notag & \qquad =  
2 \int_{-1}^1\int_{\pi+\smallasin h}^{2\pi-\smallasin h}\!\!
     \int_{-\smallasin h}^{\pi+\smallasin h}\!\!
      \int_{-\smallasin h}^{\beta_2}
\!\!\!\!\!\! \!\!    \det \betamatrix 
      d\beta_1  d\beta_2  d\beta_3 dh
\\ & \qquad
= \frac{512}{9}
\label{int:aire-arcsin} 
\\ \notag & 
 \int_{-1}^1 
      \int_{\pi+\smallasin h}^{2\pi-\smallasin h}
     \int_{-\smallasin h}^{\pi+\smallasin h}
      \int_{-\smallasin h}^{\beta_2}
     \det\betamatrix
\\ \notag  &              \hspace*{6cm}                \times
              2\sin\frac{\beta_2-\beta_1}{2}
              d\beta_1   d\beta_2d\beta_3 dh \qquad\qquad\qquad
\\ & \qquad 
= \frac{35\pi}{3}
\label{int:aireXlength-arcsin}  
\\\notag &
\int_{\frac{\pi}{2}}^{\frac{3\pi}{2}}
     \int_{\pi-\beta_3}^{\frac{\pi}{2}}
     \int_{\pi-\beta_2-2\alpha_{\rho,\kappa}}^{\pi-\beta_2+2\alpha_{\rho,\kappa}}
     \det\betamatrix  (\cos \beta_1 - \cos \beta_2)
     {d\beta_1d\beta_2d\beta_3}
\\ & \qquad =
 \tfrac{64}{9}  \cos^3  \alpha_{\rho,\kappa} \sin \alpha_{\rho,\kappa}  +\tfrac {32}{3}\cos\alpha_{\rho,\kappa} \sin\alpha_{\rho,\kappa}
+\tfrac{32}{3}\alpha_{\rho,\kappa}
\label{int:aire-projection}
\end{align}

\end{document}